\documentclass[11pt]{article}
\usepackage{amsfonts}
\usepackage{amssymb}
\newcommand{\be}{\begin{equation}}
\newcommand{\ee}{\end{equation}}
\newcommand{\ba}{\begin{array}}
\newcommand{\ea}{\end{array}}
\newcommand{\bea}{\begin{eqnarray}}
\newcommand{\eea}{\end{eqnarray}}
\newcommand{\bee}{\begin{eqnarray*}}
\newcommand{\eee}{\end{eqnarray*}}

\newtheorem{Thm}{Theorem}
\newtheorem{Lemma}{Lemma}
\newtheorem{Prop}{Proposition}
\newtheorem{Definition}{Definition}

\newtheorem{remark}{Remark}

\def\C{{\mathbb C}}

\def\R{{\mathbb R}}
\def\N{{\mathbb N}}

\def\lim{\mathop{\rm lim}}

\def\sup{\mathop{\rm sup}}

\def\aa{\alpha}
\def\dde{\delta}
\def\aad{\dot \alpha}
\def\aadd{\ddot \alpha}
\def\e{\varepsilon}
\def\l{\lambda}

\def\calH{{\cal H}}

\def\calN{{\cal N}}

\def\ld{\dot{\lambda}}
\def\bd{\dot{\beta}}

\def\gd{\dot{\gamma}}
\def\Rb{\overline{R}}

\title{Two soliton solutions to the three dimensional gravitational Hartree equation}

\author{Joachim Krieger$^{*}$, Yvan Martel$^{**}$ and Pierre Rapha\"el$^{***}$}

\date{\it $^{*}$ University of Pennsylvania, USA\\ 
$^{**}$ Universit\'e de Versailles-Saint-Quentin-en-Yvelines, France\\
$^{***}$ Universit\'e Paul Sabatier, Toulouse, France}

\begin{document}

\maketitle

\begin{abstract} We construct non dispersive two soliton solutions to the three dimensional gravitational Hartree equation whose trajectories asymptotically reproduce the nontrapped dynamics of the gravitational two body problem.
\end{abstract}


\section{Introduction}



\subsection{Setting of the problem}


We consider in this paper the three dimensional gravitational Hartree equation:
\be
\label{hartree}
 \left   \{ \begin{array}{lll}
         iu_t + \Delta u-\phi_{|u|^2}u=0,\\ 
         \Delta\phi_{|u|^2}=|u|^2 \ \ \mbox{i.e.} \ \ \phi_{|u|^2}=-\frac{1}{4\pi\|x\|}\star |u|^2,\\
         (t,x)\in \R\times \R^3, \  \ u(0,x)=u_0(x), \ \ u_0:\R^3\to \C.
         \end{array}
\right .
\ee
This system arises in Physics as an effective evolution equation in the mean field limit of many body quantum systems, see for example \cite{EY}, \cite{ES}.

It is well-known, see  \cite{Cbook} and references therein, that the Cauchy problem for (\ref{hartree}) is globally well-posed in the energy space $H^1=\{u,\nabla u\in L^2(\R^3)\}$, i.e., for $u_0\in H^1$, there exists a unique global solution $u(t)\in {\cal{C}}(\R,H^1)$ of (\ref{hartree}).

Moreover, the following quantities are conserved by the $H^1$ flow:
$$         L^2\ \ \mbox{norm}: \ \ \int|u(t,x)|^2dx=\int|u_0(x)|^2dx,$$
        $$ \mbox{Hamiltonian}:\ \ {\cal{H}}(u(t,x))=\frac{1}{2}\int|\nabla u(t,x)|^2dx-\frac{1}{4}\int |\nabla \phi_{|u|^2}(t,x)|^2dx={\cal{H}}(u_0),$$
        $$ \mbox{Momentum}: \ \ Im\left ( \int\nabla u\, \overline{u}(t,x)dx\right )=Im\left ( \int\nabla u_0\, \overline{u_0}(x)dx\right ).$$
        Recall also that the conservation of the energy and the $L^2$ norm imply an a priori bound on the kinetic energy as derived from the 
        Hardy-Littlewood-Sobolev inequality and then the Gagliardo-Nirenberg inequality:
         $$\forall u\in H^1, \ \ \int |\nabla \phi_{|u|^2}|^2
         \leq C \|u\|_{L^{\frac {12} 5}}^4
         \leq C\|\nabla u\|_{L^2}\|u\|_{L^2}^3.$$ 
This implies that (\ref{hartree}) is subcritical so that all $H^1$ solutions to (\ref{hartree}) are global and bounded in $H^1$. 
In contrast, note that in dimension four, the Hartree equation with potential $\frac 1{\|x\|^2}$ is $L^2$ critical and a singularity formation is known to possibly occur, see \cite{KLR} for recent progress on this problem.

Equation (\ref{hartree}) possesses a large group of symmetries: if $u(t,x)$ is a solution to (\ref{hartree}), then for all $(t_0,\aa_0,\beta_0,\gamma_0,\lambda_0)\in \R\times\R^N\times\R^N\times \R\times \R_*^+$, so is 
\be
\label{symmetry}
v(t,x)=\lambda_0^{2}u(\lambda_0^2 t+t_0,\lambda_0 x+\aa_0-\beta_0 t)e^{i\frac{\beta_0}{2}\cdot(x-\frac{\beta_0}{2} t)}e^{i\gamma_0}.
\ee  

Special solutions are expected to play a fundamental role for the description of the dynamics of (\ref{hartree}). They are the so-called solitary waves, of the form $u(t,x)=e^{it}W(x)$ where $W(x)$ solves
\be
\label{eqsoliton}
\Delta W - \phi_{|W|^2}W= W.
\ee
We shall denote by $Q$ the so-called ground state solution to (\ref{eqsoliton}), which is defined as the unique radially symmetric nonnegative solution to (\ref{eqsoliton}).  The existence and {\it uniqueness} of the ground state has been proved using variational and ODE techniques by Lieb \cite{Lieb}.

The group of symmetries (\ref{symmetry}) then generates an eight parameter family of ground state solitary waves: $$Q_{\aa,\beta,\gamma,\lambda}(t,x)=\lambda^{2}Q(\lambda x+\aa-\beta t)e^{i\frac{\beta}{2}\cdot(x-\frac{\beta}{2} t)}e^{i\lambda^2 t}e^{i\gamma},  \ \ (\aa,\beta, \lambda,\gamma)\in \R^3\times\R^3\times \R \times \R_+^* .$$ Observe that the center of mass of such a solution  is evolving according to the free Galilean motion with constant speed $\beta$.
Recall that a variational characterization of the ground state $Q$ has been derived by Lions \cite{PLL1} using the concentration-compactness technique. Given $M>0$, the minimization problem $$\inf_{\|u\|_{L^2}=M}{\cal{H}}(u)$$ is attained at 
$$u(x)=\lambda^2(M)Q(\lambda(M) x+\aa)e^{i\gamma}, \  \ (\aa,\gamma)\in \R^3\times \R, \ \ \lambda(M)=\frac{M^2}{\|Q\|_{L^2}^2}.$$
Moreover, every minimizing sequence is relatively compact in the energy space up to phase and translation shifts. Following Cazenave and Lions \cite{CL} and using Lieb's  uniqueness result \cite{Lieb},  this automatically implies the orbital stability of the ground state solitary wave.

\medskip

The question of the long time dynamics of the Hartree equation is widely open. However, in other settings like the one dimensional (gKdV) equation or the Schr\"odinger equation (NLS) with  power nonlinearity 
\be
\label{nls}
  \left   \{ \begin{array}{ll}
         iu_t=-\Delta u-|u|^{p-1}u, \ \ (t,x)\in \R\times \R^N\\
         u(0,x)=u_0(x), \ \ u_0:\R^N\to \C, \ \ p\leq 1+\frac{4}{N},
         \end{array}
\right .
\ee
multisolitary wave solutions are conjectured to be the building blocks for the description of the long time dynamics. Roughly speaking, as time goes to $+\infty$, a generic solution should split asymptotically into a sum of solitary waves which move away from each other and a radiative term which disperses. There is a very extensive literature on this subject, we refer for example to \cite{BP}, \cite{RSS}, \cite{MMTkdv}, \cite{MMTnls}, \cite{Cote}, \cite{TY},  \cite{Tao} (and references therein).

\medskip

It seems that one of the first problem to address for the understanding of the long time dynamics is the existence of {\it non dispersive} multisolitary waves. Apart from completely integrable systems like the (KdV) equation  or the cubic one dimensional (NLS) where explicit multisolitary waves can be exhibited (see e.g. \cite{Z}), the existence of such objects goes back to Merle \cite{Mcrit} for critical (NLS) problems, and has later been extended by Martel \cite{Mkdv} to subcritical (KdV) equations and by Martel and Merle \cite{MMnls} for general subcritical nonlinear (NLS) problems. See also Buslaev and Perelman \cite{BP} and Rodnianski, Soffer and Schlag  \cite{RSS} for  different approaches. 

Let us observe that for a system like the subcritical (NLS) (\ref{nls}), the sum of two solitary waves moving away from each other at constant speed is a solution to (\ref{nls}) up to an exponentially small in time correction. This suggests that there is essentially no interaction between the two ground states, and that nondispersive multisolitary wave solutions can be constructed with asympotic trajectories given by the free Galilean motion along any chosen non parallel lines, \cite{MMnls}.


\subsection{Statement of the result}


This paper concerns the construction of two soliton solutions of the Hartree equation.
Note that for the Hartree equation, in contrast with the (NLS) or (gKdV) case, even though the ground state is exponentially decreasing in space, the nonlinearity is {\it long range} from the slow decay of the gravitational field $$\phi_{|Q|^2}(r)\sim\frac{C}{r}\ \ \mbox{as} \ \ r\to+\infty.$$
The consequence is  a strong coupling between the two solitons. In this paper, we claim that this constraint implies a non trivial dynamics for the centers of mass of the solitons, which -- as one should expect from physical grounds -- asymptotically converges to the dynamical system of the two body problem in Newtonian gravity. 

Now, we recall the well-known definition of the two body problem:

\begin{Definition}[Dynamics of the two body problem]\quad 
A trajectory of the two body problem with masses $(\lambda_1,\lambda_2)\in (\R^*_+)^2$ is a solution $(\aa_1(t),\aa_2(t),\beta_1(t),\beta_2(t))$ to 
\be
\label{twobodyhyp}
 \left   \{ \begin{array}{ll}
         \aad_1=2\beta_1, \ \ \aad_2=2\beta_2,\ \ \aa=\aa_2-\aa_1, \\
         \bd_1=\frac{\|Q\|_{L^2}^2}{4\pi\lambda_2}\frac{\aa}{\|\aa\|^3}, \ \ \bd_2=-\frac{\|Q\|_{L^2}^2}{4\pi\lambda_1}\frac{\aa}{\|\aa\|^3}, \ \ \beta=\beta_2-\beta_1.
         \end{array}
\right .
\ee
The center of mass evolves according to the free Galilean motion, i.e. $\lambda_2\aadd_1+\lambda_1\aadd_2=0$ and $\aa(t)$ evolves in a fixed plane. Moreover, the dynamical system admits the conserved Hamiltonian:
$$E_0=\|\beta\|^2-\frac{\|Q\|_{L^2}^2}{4\pi}\left(\frac{1}{\lambda_1}+\frac{1}{\lambda_2}\right)\frac{1}{\|\aa\|}.$$ 
The dynamics are classified as follows:

(i) Hyperbolic trajectory: If $E_0>0$ then $\aa(t)$ describes a hyperbola with \be
\label{limiteinfinie}
\lim_{t\to +\infty}\frac{\|\aa(t)\|}{t}=2\sqrt{E_0}. 
\ee

(ii) Parabolic trajectory:  If $E_0=0$ then $\aa(t)$ describes a parabola with
\be
\label{limiteinfinieparabole}
\lim_{t\to +\infty
}\frac{\|\aa(t)\|}{t^{\frac{2}{3}}}=C>0.
\ee

(iii) Elliptic trajectory: If $E_0<0$ then the dynamic is periodic in time and $\aa(t)$ describes an ellipse.
\end{Definition}

In the rest of this paper, we fix $(\lambda_1^{\infty},\lambda_2^{\infty})\in (\R^*_+)^2$, and we consider a given solution $(\aa^{\infty}_1(t), \aa^{\infty}_2(t), \beta^{\infty}_1(t),\beta^{\infty}_2(t))$ of the two body problem with masses  $(\lambda_1^{\infty},\lambda_2^{\infty})$.
We assume throughout the paper that the center of mass is fixed at the origin and that  the trajectory lies on the plane $\{x_{(3)}=0\}$ where $x=(x_{(1)},x_{(2)},x_{(3)})$, i.e.
\be
\label{assumption}
\forall t\geq 0, \ \ \lambda_2^{\infty}\aa_1^{\infty}(t)+\lambda_1^{\infty}\aa_2^{\infty}(t)=0  \quad \hbox{and} \quad Ê(\alpha_{j}^{\infty})_{(3)}(t)=0,
\ee
where $\alpha_j(t)=((\alpha_j)_{(1)}(t), (\alpha_j)_{(2)}(t),(\alpha_j)_{(3)}(t)).$ Note that these assumptions do not restrict generality by standard use of translation, rotation and Galilean invariances.

\medskip

Our main result concerns the existence of {\it non dispersive} two soliton solutions which asymptotically reproduce the {\it non trapped dynamics} of the two body problem in the following two cases:
\begin{enumerate}
\item Hyperbolic case without restriction on the masses;
\item Parabolic case with equal masses.
\end{enumerate}

\begin{Thm}[Existence of a two soliton with a non trapped trajectory]
\label{thmhyperbolic}\quad \\
Let $(\lambda^{\infty}_1,\lambda^{\infty}_2)\in (\R^*_+)^2$ and let $(\aa^{\infty}_1(t), \aa^{\infty}_2(t), \beta^{\infty}_1(t),\beta^{\infty}_2(t))$ be a solution to the two body problem with masses $(\lambda^{\infty}_1,\lambda^{\infty}_2)$
such that  (\ref{assumption}) holds.
\begin{enumerate}
\item Hyperbolic case.
Assume that  $(\aa^{\infty}_1(t), \aa^{\infty}_2(t), \beta^{\infty}_1(t),\beta^{\infty}_2(t))$ is hyperbolic $(E_0>0)$. Then, there exists an $H^1$ solution of (\ref{hartree}) and $\gamma_1(t),\gamma_2(t)$
such that
\be\label{thm1-1}
\lim_{t\to +\infty}\left\|u(t,x)- \sum_{j=1}^2 \frac{1}{(\lambda_j^\infty)^2}Q\left(\frac{x-\aa_j^\infty(t)}{{\lambda_j^\infty}}\right)e^{-i\gamma_j(t)+i\beta_j^\infty (t)\cdot x}\right\|_{H^1}= 0.
\ee
\item Parabolic case. Assume that $(\aa^{\infty}_1(t), \aa^{\infty}_2(t), \beta^{\infty}_1(t),\beta^{\infty}_2(t))$  is  parabolic $(E_0=0)$ and $\lambda_1^\infty=\lambda_2^\infty=\lambda^\infty$.
Then, there exist an $H^1$ solution of (\ref{hartree}) and $\alpha_1(t),$ $\alpha_2(t)$, $\gamma_1(t),$ $\gamma_2(t)$ such that
$$
\lim_{t\to +\infty}\left\|u(t,x)- \sum_{j=1}^2 \frac{1}{(\lambda^\infty)^2}Q\left(\frac{x-\aa_j(t)}{{\lambda^\infty}}\right)e^{-i\gamma_j(t)+i\beta_j^\infty (t)\cdot x}\right\|_{H^1}= 0,
$$
$$
\lim_{t\to +\infty} \sum_{k=1}^2\left|\frac{(\aa_j^{\infty})_{(k)}(t)}{(\aa_j)_{({k})}(t)}-1\right|+|(\alpha_j)_{(3)}(t)| = 0,  \ \ j=1,2.
$$
\end{enumerate}
\end{Thm}

\noindent{\it Comments on the result:}
\smallskip

{\it 1. Constraints on the asymptotic trajectory.} Note that for local nonlinearities with weak interactions, like in the case of power nonlinearity (NLS) problems, the asymptotic dynamics of the centers of mass can be prescribed along any prescribed lines, see \cite{MMnls}. In contrast, in Theorem \ref{thmhyperbolic}, the asymptotic hyperbolic trajectory has to be along two {\it coplanar} lines. We expect this situation to be the only possibility.

\smallskip

{\it 2. Parabolic case.} The behavior displayed by the solution constructed in Theorem 1 in the parabolic case does not correspond asymptotically to an almost free Galilean motion of each body, i.e. along two coplanar lines, which is for one soliton the natural motion from the Galilean symmetry. Thus, the result of Theorem 1 in the parabolic case 
displays a new nonlinear (critical) regime where the trajectories of the two solitary waves do not correspond to asymptotically free solitary waves. 

Note also from the proof that we do not know whether $\|\aa_j(t)-\aa_j^{\infty}(t)\|\to 0$ as $t\to +\infty$ in the parabolic case. This question is related to the computation of a constant, which does not seem to be explicit. However, if we set ${\cal{P}}^\infty=\{\aa_1^{\infty}(t), \aa_2^{\infty}(t), \ t\geq 0\}$, then  $$\lim_{t\to +\infty} \mbox{dist}((\aa_1(t),\aa_2(t)),{\cal{P}}^{\infty})=0,$$
which means that $(\aa_1(t),\aa_2(t)))$ asymptotically describes the same set as $(\aa_1^{\infty}(t), \aa_2^{\infty}(t))$ with different speed. See Appendix A.2.

\smallskip

{\it 3. Connection to existing literature.} The fact that the Hartree problem displays solutions which reproduce in some sense the point particle newtonian interaction is not a surprise. This has been in particular observed by Frohlich, Yau, Tsai \cite{FTY} where effective modulation equations leading to the Newtonian Galilean motion for the center of mass are derived in certain asymptotic regimes  -- but not for all times $t\to +\infty.$ In \cite{SZ}, Gang and Sigal describe the long time behavior of the (NLS) equations in a trapping very well localized potential and show how the Newtonian law of motion is here again driving the center of mass of the ground state part of the solution.

\smallskip

{\it 4. Uniqueness.} For the hyperbolic case, one may ask whether a solution satisfying (\ref{thm1-1})
is unique. By the method of this paper, one can expect a weaker statement, i.e. uniqueness in a class  smaller than (\ref{thm1-1}).
Uniqueness in general is an open problem, certainly related to the -delicate- stability problem, see next comment. See \cite{Mkdv} for a general uniqueness proof in the (gKdV) case. 

\smallskip

{\it 5. Stability.} The main open problem after this work is the stability of the two soliton solution as derived for some (NLS) and (gKdV) type problems, \cite{RSS}, \cite{MMTkdv}, \cite{MMTnls}. We expect the instability of the parabolic trajectory and the stability of the hyperbolic one. Another important open problem concerns the non existence of trapped elliptic type trajectories which most likely are  destroyed by the dispersive effects of the flow.

\medskip

The proof of Theorem \ref{thmhyperbolic} relies on a refinement of the techniques developed by Martel and Merle \cite{Mkdv}, \cite{MMnls}, and proceeds into two steps.
\begin{itemize}
\item Construction of an approximate solution to all orders: The first step is to display a general framework to compute an approximate solution to the two soliton equation with an arbitrarily high order error $\frac{1}{\|\aa\|^N}$ where $\|\aa\|$ is the distance between the center of masses. Here our analysis relies on a {\it finite dimensional reduction} of the problem. This kind of procedure is reminiscent to the computation of formal ansatz and solutions in various nonlinear dispersive settings including for example the description of singularity formations, \cite{SS}, \cite{MR2}, \cite{Kr-Sch-Tat},  the description of the long time dynamics of the (NLS) near a ground state, \cite{SZ}, or the interaction of two solitons \cite{MMinter}. We propose here a systematic way of producing such approximate solutions. The main advantage is that the approximate solution is now good enough to somehow reduce the problem to an almost {\it short range} problem.
\item Construction of the nondispersive two soliton: The second step is to build the exact two soliton solution by solving the problem for the small remainder backwards from infinity with a procedure which is reminiscent from the construction of nonlinear wave operators -see for example \cite{BW} for a similar strategy-. The control of the remainder relies on a space localization of the conservation laws and {\it energy estimates} on the linearized flow close to a solitary wave. Such estimates are a consequence of the variational characterization of $Q$ as first exhibited by Weinstein \cite{W1} in the local (NLS) setting, and some elliptic non degeneracy properties proved by Lenzman \cite{Lenzman} for the Hartree problem. While this strategy has the advantage of not requiring any dispersive estimate, the localization of conservation laws creates large errors in our setting due to the long range structure of the problem and the only polynomial distance between the two solitary waves. In fact, it turns out that the closure of these estimates is critical here with respect to the $\frac{1}{r}$ decay of the gravitational field, in particular in the parabolic case where we can only treat the case of an asymptotic {\it symmetric} two body problem. Eventually, the success of our strategy relies in a crucial way onto the the construction of an arbitrary large order approximate solution from the first step, see Lemma \ref{energyestimate} and Remark \ref{rkkey}.
\end{itemize}

We expect this strategy to be quite robust and to apply to a very large class of situations and problems.

\medskip

\noindent{\bf Acknowledgements.} J. K. is supported by NSF-Grant DMS-757278 and a Sloan Fellowship. Y.M is supported by ANR Projet Blanc OndeNonLin. P.R is supported by ANR jeunes chercheurs SWAP. Part of this work was done while J.K was visiting the Institut de Math\'ematiques de l'Universit\'e Paul Sabatier, Toulouse, which he would like to thank for its kind hospitality.

\section{Construction of the approximate solution}

The purpose of this section is to construct a formal approximate two soliton solution of (\ref{hartree}) at any order of $\frac 1{\|\aa\|}$, see Proposition \ref{propinduction} and  Proposition \ref{constructionansatz}. Our strategy relies on a finite dimensional reduction of the problem which is connected to the strategies developed for example in \cite{SS}, \cite{MR2}, \cite{SZ}, \cite{MMinter}, \cite{Kr-Sch-Tat}.



\subsection{Derivation of the coupled equations}


Let us start with computing the evolution equation in renormalized variables. Let 
$$u(t,x)=\frac{1}{\lambda^2(t)}v\left(t,\frac{x-\aa(t)}{\lambda(t)}\right)e^{-i\gamma(t)}e^{i\beta(t)\cdot x},\ \ v=v(t,y).$$ Then 
\bea
\label{changevariables}
& & i\partial_t u+\Delta u-\phi_{|u|^2}u= \frac{1}{\lambda^4}\Big[i\lambda^2\partial_tv+\Delta v-v-\phi_{|v|^2}v\\
\nonumber&  &  -i\lambda\ld\Lambda v -\lambda^3(\bd\cdot y)v-i\lambda(\aad-2\beta)\cdot\nabla v\\
\nonumber & &  + \lambda ^2(\gd+\frac{1}{\lambda^2}-\|\beta\|^2-\bd\cdot \aa)v\Big]\left(t,\frac{x-\aa(t)}{\lambda(t)}\right)e^{-i\gamma(t)}e^{i\beta(t)\cdot x},
\eea
where $\Lambda$ denotes the
differential operator $$\Lambda v=2v+y\cdot\nabla v.$$ 
We compute an ansatz of a two soliton solution as $t\to +\infty$ by letting 
\bee
u(t,x) & = & u_1(t,x)+u_2(t,x)\\
& = &\frac{1}{\lambda^2_1(t)}v_1\left(t,y_1\right)e^{-i\gamma_1(t)}e^{i\beta_1(t)\cdot x}+\frac{1}{\lambda^2_2(t)}v_2\left(t,y_2\right)e^{-i\gamma_2(t)}e^{i\beta_2(t)\cdot x}
\eee where
\be
\label{rescaledcoord}
y_j=\frac{x-\aa_j(t)}{\lambda_j(t)}, \ \ j=1,2.
\ee 
The nonlinear term is $$\phi_{|u|^2}u=\phi_{|u_1+u_2|^2}(u_1+u_2).$$ By construction, $u_1$ and $u_2$ will have disjoint supports up to exponentially small corrections. Hence we will not consider the crossed term for now (see Proposition \ref{constructionansatz} for details). By rescaling, $$\phi_{|u_j|^2}(x)=\frac{1}{\lambda_j^2}\phi_{|v_j|^2}(y_j), \ \ j=1,2$$ and so
\bee
&& \phi_{|u|^2}u  =  (\phi_{|u_1|^2}+\phi_{|u_2|^2})u_1+(\phi_{|u_1|^2}+\phi_{|u_2|^2})u_2+O(e^{-\gamma t})\\
& = & \frac{1}{\lambda_1^4}\left[\phi_{|v_1|^2}+\left(\frac{\lambda_1}{\lambda_2}\right)^2\phi_{|v_2|^2}\left(\frac{\lambda_1}{\lambda_2}y_1-\frac{\aa_2-\aa_1}{\lambda_2}\right)\right]v_1(t,y_1)e^{-i\gamma_1(t)}e^{i\beta_1(t)\cdot x}\\
& + & \frac{1}{\lambda_2^4}\left[\phi_{|v_2|^2}+\left(\frac{\lambda_2}{\lambda_1}\right)^2\phi_{|v_1|^2}\left(\frac{\lambda_2}{\lambda_1}y_2+\frac{\aa_2-\aa_1}{\lambda_1}\right)\right]v_2(t,y_2)e^{-i\gamma_2(t)}e^{i\beta_2(t)\cdot x}+O(e^{-\gamma t}).
\eee
Thus, we rewrite the equation $$ i\partial_tu+\Delta u-\phi_{|u|^2}u= \sum_{j=1}^2 \frac{1}{\lambda^4_j}Eq_j(t,x)e^{-i\gamma_j(t)}e^{i\beta_j(t)\cdot x}$$ with 
\bee
Eq_j(t,x) & = &  i\lambda_j^2\partial_tv_j+\Delta v_j-v_j-\left[\phi_{|v_j|^2}+\phi_{|v_{j+1}|^2}\left(\frac{\lambda_j}{\lambda_{j+1}}y_{j}+(-1)^j\frac{\aa_2-\aa_1}{\lambda_{j+1}}\right) \right]v_j\\
\nonumber & - & i\lambda_j\ld_j\Lambda v_j -\lambda_j^3(\bd_j\cdot y_j) v_j-i\lambda_j(\aad_j-2\beta_j)\cdot\nabla v_j\\
\nonumber& + & \lambda_j ^2(\gd_j+\frac{1}{\lambda_j^2}-\|\beta_j\|^2-\bd_j\cdot \aa_j)v_j
\eee
and the convention of mod(2) summation: $$j+1=1\ \ \mbox{for}  \ \ j=2.$$

\subsection{Dipole expansion of the gravitational field}\label{sec:1.2}


Let
\be
\label{defxbeta}
\aa(t)=\aa_2(t)-\aa_1(t), \ \ \beta(t)=\beta_2(t)-\beta_1(t).
\ee
We expand the gravitational field created by the soliton $u_{j+1}$ on $u_j$:
$$
  \left(\frac{\lambda_j}{\lambda_{j+1}}\right)^2\phi_{|v_{j+1}|^2}\left(\frac{\lambda_j}{\lambda_{j+1}}y_{j}+(-1)^j\frac{\aa}{\lambda_{j+1}}\right) = -\frac{\lambda_j^2}{4\pi\lambda_{j+1}}\int_{\R^3}\frac{|v_{j+1}(\xi)|^2}{\|\aa-(-1)^j(\lambda_{j+1}\xi-\lambda_jy_j)\|}d\xi.
 $$
We proceed to the formal dipole expansion of the gravitational potential which main contribution corresponds to the region  $\|y_j\|+\|\xi\|\ll\|\aa\|$. Such a formal argument will be justified in the proof of Proposition \ref{constructionansatz} below for functions $v_j$ with sufficiently fast decay at infinity in space. 

Let 
\be\label{pretaylor}
\frac{1}{\|\aa-\zeta\|}= \frac 1 {\|\aa\|} \left(1- 2 \frac {\zeta\cdot \aa}{\|\aa\|^2} 
+ \frac {\|\zeta\|^2}{\|\aa\|^2}\right)^{-\frac{1}{2}}
\ee and use a Taylor expansion:
\be
\label{expansionofthefield}
\frac{1}{\|\aa-\zeta\|}\sum_{k=1}^N F_k(\aa,\zeta)+O\left(\frac{\|\zeta\|^{N+1}}{\|\aa\|^{N+1}}\right).
\ee
Explicit computations yield for the first terms: $$F_1(\aa,\zeta)=\frac{1}{\|\aa\|}, \ \ F_2(\aa, \zeta)=\frac{\aa \cdot \zeta}{\|\aa \|^3} , \ \ F_3(\aa,\zeta)=\frac{1}{\|\aa\|^3}\left[\frac{3}{2}\left(\frac{\aa}{\|\aa\|}\cdot\zeta\right)^2-\frac{1}{2}\|\zeta\|^2\right].$$ Observe in particular that $F_k$ is homogeneous of degree $-k$ in $\aa$.\\
We thus introduce the approximated correction field of order $N$: 
\be
\label{approxfield}
\phi_{|v_{j+1}|^2}^{N,app}(y_j)=\sum_{k=1}^N\phi_{|v_{j+1}|^2}^{N,k,app}(y_j)
\ee
with 
\be
\label{fieldnk}
\phi_{|v_{j+1}|^2}^{N,k,app}(y_j)=-\frac{\lambda_{j}^2}{4\pi\lambda_{j+1}} \int_{\R^3}|v_{j+1}(\xi)|^2F_k((-1)^j \aa,\lambda_{j+1}\xi-\lambda_j y_j)d\xi.
\ee


\subsection{Finite dimensional reduction of the dynamics}


From the previous formal calculations, we aim at constructing an approximate solution to the system: $Eq^{(N)}_j=0$ for $j=1,2$, where
\bee
Eq^{(N)}_j
& = &  i\lambda_j^2\partial_tv_j+\Delta v_j-v_j-\left[\phi_{|v_j|^2}+\phi_{|v_{j+1}|^2}^{N,app}\right]v_j\\
\nonumber & - & i\lambda_j\ld_j\Lambda v_j -\lambda_j^3(\bd_j\cdot y_j) v_j-i\lambda_j(\aad_j-2\beta_j)\cdot\nabla v_j\\
\nonumber& + & \lambda_j ^2(\gd_j+\frac{1}{\lambda_j^2}-\|\beta_j\|^2-\bd_j\cdot \aa_j)v_j.
\eee
We now proceed to a finite dimensional reduction of this problem and look for a solution of the form $$v^{(N)}_j(t,y_j)=\left[V_j^{(N)}\right]_{(\aa(t),\beta(t),\lambda_1(t),\lambda_2(t))}(y_j).$$  Here $V_j^{(N)}$ is stationary and hence $v_j^{(N)}$ depends on time only through the modulation parameters $t\to (\aa(t),\beta(t),\lambda_1(t),\lambda_2(t))$ to be chosen: 
$$\partial_t v_j^{(N)}= \sum_{k=1}^2\frac{\partial V_j^{(N)}}{\partial\lambda_k}\ld_k+\frac{\partial V_j^{(N)}}{\partial \aa}\cdot\aad+\frac{\partial V_j^{(N)}}{\partial\beta}\cdot\bd.$$ We obtain 
\bea
\label{approximationone}
\nonumber Eq^{(N)}_j&= &  \Delta V^{(N)}_j-V^{(N)}_j-\left[\phi_{|V^{(N)}_j|^2}+\phi_{|V_{j+1}^{(N)}|^2}^{N,app}\right]V^{(N)}_j \\ &-&    i\l_jM_j^{(N)}\Lambda V^{(N)}_j -\lambda_j^3B_j^{(N)}\cdot y_j V^{(N)}_j\\
\nonumber &+& i \lambda_j^2 \left[\sum_{k=1}^2M_k^{(N)}\frac{\partial V_j^{(N)}}{\partial\lambda_k}+\frac{\partial V_j^{(N)}}{\partial \aa}\cdot 2\beta + \frac{\partial V_j^{(N)}}{\partial \beta}\cdot (B_2^{(N)}-B_1^{(N)})\right] +S_j^{(N)},
\eea
where $S_j^{(N)}$ encodes the finite dimensional system of the geometrical parameters:
\bea
\label{defsnj}
\nonumber S_j^{(N)} & = & - i \lambda_j(\aad_j-2\beta_j)\cdot\nabla V^{(N)}_j-  \lambda_j^3(\bd_j-B_j^{(N)})\cdot y_jV_j^{(N)} \\
\nonumber & - &i\lambda_j(\ld_j-M_j^{(N)})\Lambda V_j^{(N)}+  \lambda_j ^2(\gd_j+\frac{1}{\lambda_j^2}-\|\beta_j\|^2-\bd_j\cdot \aa_j)V^{(N)}_j\\
\nonumber & + &i \lambda_j^2\sum_{k=1}^2\left[(\ld_k-M_k^{(N)})\frac{\partial V_j^{(N)}}{\partial\lambda_j}+(-1)^{k}(\aad_{k}-2\beta_{k})\cdot\frac{\partial V_j^{(N)}}{\partial \aa}\right .\\
& +& \left . (-1)^k(\bd_k-B_k^{(N)})\cdot \frac{\partial V_j^{(N)}}{\partial\beta}\right]
\eea
Now, we build an approximate solution using a series expansion of $V_j^{(N)}$ in $(\aa,\beta,\lambda_1,\lambda_2)$ with a prescribed structure. We need the following definition. 

\begin{Definition}[Admissible functions]
\label{admissible}
 Let $n\in \N$. \\
 (i) We define ${\cal{S}}_{n}$  the class of complex valued functions $\sigma(\aa,\beta,\lambda_1,\lambda_2):(\R^3\setminus\{0\})\times(\R^3\setminus\{0\})\times\R^*_+\times\R^*_+\to \C$ of the form:
 $$\sigma(\aa,\beta,\lambda_1,\lambda_2)=\frac{1}{\lambda_1^{n_1}\lambda_{2}^{n_2}}\sum_{l_1,l_2,k_1,k_2,k_3=0}^{n_3}\lambda_1^{l_1}\lambda^{l_2}_{2}\beta_{(1)}^{k_1}\beta_{{(2)}}^{k_2}\beta_{(3)}^{k_3}f_{l_1,l_2,k_1,k_2,k_3}(\aa),\ \ n_1,n_2,n_3 \in \N,$$ where 
 $\beta=(\beta_{(1)},\beta_{(2)},\beta_{(3)})$ and $f_{l_1,l_2,k_1,k_2,k_3}$ is a complex valued ${\cal{C}}^{\infty}$ function of $\aa$ away from $\aa=0$ and homogeneous of degree $-n$ in $\aa$. Then, we set $$deg(\sigma)=n.$$
 (iii) We say that a function $\psi: (\aa,\beta,\lambda_1,\lambda_2,x):(\R^3\setminus\{0\})\times(\R^3\setminus\{0\})\times\R^*_+\times\R^*_+\times\R^3\to \C$ is admissible of degree $n$  if it admits a finite expansion of the form $$\psi(\aa,\beta,\lambda_1,\lambda_2,x)= \sum_{m=0}^{n_4}\sigma_{m}(\aa,\beta,\lambda_1,\lambda_2)\tau_{m}(x), \ \ n_4\in \N$$ where $\sigma_{m}\in {\cal{S}}_{n}$ and $\tau_{m}$ is well-localized in $x$: 
\be
\label{expdecay}
\forall k\geq 0, \forall x\in \R^3, \ \ \left\|\nabla^k \tau_{m}(x)\right\|\leq e^{-C_{m,k}\|x\|}.
\ee
Then, we set $$deg(\psi)=n.$$
(iv) We say that a function $\psi: (\aa,\beta,\lambda_1,\lambda_2,x):(\R^3\setminus\{0\})\times(\R^3\setminus\{0\})\times\R^*_+\times\R^*_+\times\R^3\to \C$ is admissible of degree $\geq N$ if it admits a finite expansion of the form $$\psi=\sum_{n=N}^{n_5}\psi_n,\ \ n_5\geq N$$ with $\psi_n$ admissible of degree $n$. We say that $deg(\psi)\geq N$ with abuse of notation.
\end{Definition}

We claim the following technical facts.

\begin{Lemma}  
\label{stabadmissible}
Let $(n_1,n_2)\in \N^2$. Let $\psi_1,\psi_2$ be admissible of degree respectively $n_1,n_2$. Then:

(i) $\psi_1\psi_2$ is admissible of degree $n_1+n_2$.

(ii) $\phi_{\psi_1}\psi_2$ is admissible of degree $n_1+n_2$.

(iii) $\forall N\geq 1$, $\phi^{N,app}_{\psi_1}\psi_2$ is admissible of degree $\geq 1+n_1+n_2$.

(iv) $\forall 1\leq k\leq N$, $\phi^{N,k,app}_{\psi_1}\psi_2$ is admissible of degree $k+n_1+n_2$.
\end{Lemma}

{\bf Proof.}
Properties (i)-(ii) are clear by the definition of admissible function of degree $n$ 
(see Definition \ref{admissible}) and the properties of homogeneous functions.
Property (iii) is a consequence of (iv). 

Finally the proof of (iv) relies
onto the fact that $F_k((-1)^j \aa , \lambda_{j+1} \xi - \lambda_{j} y_j)$ appearing in (\ref{fieldnk}) is a polynomial in $(\lambda_1,\lambda_2)$ which coefficients are homogeneous of order $-k$ in $\aa$ and $k$ in $\zeta$. Indeed, $F_k(\aa,\zeta)$ is the $k$-th term to the Taylor expansion (\ref{expansionofthefield}) of  (\ref{pretaylor}).

\medskip

Then, we claim the following Proposition. 

\begin{Prop}[Expansion of the approximate solution]
\label{propinduction}
Let $\lambda_1, \lambda_2>0$, $(\aa,\beta)\in (\R^3/\{0\})^2$ and $N\geq 1$. Let $j=1,2$ with the convention $j+1=1$ for $j=2$. We can find an asymptotic expansion 
\be
\label{computV}
V^{(N)}_j(y_j)=Q(y_j)+\sum_{n=1}^NT^{(n)}_{j}(\aa, \beta, y_j,\lambda_1,\lambda_2), 
\ee 
\be
\label{computmb}
M_j^{(N)}(\aa,\beta,\lambda_1,\lambda_2)=\sum_{n=2}^Nm_j^{(n)}(\aa,\beta,\lambda_1,\lambda_2),  \ \ B_j^{(N)}(\aa,\beta,\lambda_1,\lambda_2)=\sum_{n=2}^N b_j^{(n)}(\aa,\beta,\lambda_1,\lambda_2),\ee 
such that the following holds true.\\
(i) $T^{(n)}_{j}$ is admissible of degree $n$.\\
(ii) $m_j^{(n)}\in {\cal{S}}_{n}$, $b_j^{(n)}\in{\cal{S}}_{n}$.\\
(ii) Approximate solution of order $N+1$:  Let
\bea
\label{defeqtilden}
\nonumber\tilde{Eq}^{(N)}_j &= &  \Delta V^{(N)}_j-V^{(N)}_j-\left[\phi_{|V_j^{(N)}|^2}+\phi_{|V^{(N)}_{j+1}|^2}^{N,app}\right]V^{(N)}_j\\
\label{eq:22} &  -&   i\l_jM_j^{(N)}\Lambda V^{(N)}_j -\lambda_j^3B_j^{(N)}\cdot y_j V^{(N)}_j\\ \nonumber &  + & i\left[\sum_{k=1}^2M_k^{(N)}\frac{\partial V_j^{(N)}}{\partial\lambda_k}+\frac{\partial V_j^{(N)}}{\partial \aa}\cdot 2\beta +\frac{\partial V_j^{(N)}}{\partial \beta}\cdot (B_2^{(N)}-B_1^{(N)})\right]
\eea
Then $\tilde{Eq}^{(N)}_j$ is admissible of degree $\geq N+1$.
\end{Prop}

The proof of Proposition \ref{propinduction} proceeds by induction on $N$. We start with computing explicitely the first two  terms in powers of $\|\alpha\|^{-1}$.


\subsection{Computation of the first terms}


First, we recall the following standard facts on the linearized operator around $Q$
$$
L_+ f = - \Delta f + f + \phi_{Q^2} f + 2 \phi_{(Qf)} Q,\qquad
L_- f = - \Delta f + f + \phi_{Q^2} f .
$$
which are consequences of the nondegeneracy of the kernel of $L_+,L_-$ as exhibited in \cite{Lenzman}.

\begin{Lemma}[Inversion of $L_\pm$]
\label{lemmainversion}
Let $n\in \N$ and $f$ be real-valued, admissible of degree $n$.\\
(i) If $(f,\nabla Q)=0$, then there exists a solution to $L_+u=f$ which is real-valued and admissible of degree $n$. Moreover, if $f$ is radially symmetric then $u$ can be chosen radially symmetric. \\
(ii) If $(f,Q)=0$,  then there exists a solution to $L_-u=f$ which is real-valued and admissible of degree $n$. Moreover, if $f$ is radially symmetric then so is $u$.
\end{Lemma}  

\textbf{Proof.}
(i) Introduce the Hilbert space ${\mathcal{H}}$ of all real-valued functions $f\in L^{2}(\R^{3})$ which satisfy the orthogonality relation 
\[
(f,\, \nabla Q)=0.
\]
Then note that $L_{+}$ maps ${\mathcal{H}}\cap D$ into ${\mathcal{H}}$, $D$ being the domain of $L_{+}$ in $L^{2}(\R^{3})$. It is proved in \cite{Lenzman} that the kernel of $L_{+}$ consist precisely of the components of $\nabla Q$. Hence for $v\in{\mathcal{H}}\cap D$,  we have $L_{+}v=0$ iff $v=0$. Now consider the equation 
\[
L_{+}u=f
\]
for given $f\in {\mathcal{H}}$. Writing 
\[
u=(-\Delta+I)^{-1}v,
\]
we can write the above equation as the following problem:
\[
v+Kv=f,\ \ K =K_1 \circ (-\Delta+I)^{-1}, \ \ K_1 v= [\phi_{Q^{2}}v+2\phi_{(Qv)}Q].
\]
The operator $K: {\mathcal{H}}\rightarrow{\mathcal{H}}$ and is compact. Thus, its adjoint  $K^{*}: {\mathcal{H}}\rightarrow{\mathcal{H}}$,  is also compact. Moreover, from the information on the kernel of $L_+$, $I+K^{*}$ has a trivial kernel: indeed, denoting the projection in $L^{2}(\R^{3})$ onto $\big(span(\partial_{x_{1}}Q, \partial_{x_{2}}Q, \partial_{x_{3}}Q)\big)^{\perp}$ by $\Pi$, we have 
\[
K^{*}=\Pi\circ(-\Delta+I)^{-1}\circ K_{1}
\]
Hence the assertion $(I+K^{*})v=0$ for some $v\in  {\mathcal{H}}$ implies 
\[
v+\left((-\Delta+I)^{-1}\circ K_{1}\right) v=a\cdot\nabla Q
\]
for suitable $a\in \R^{3}$. Standard elliptic theory implies $v\in C^{\infty}(\R^{3})$, and then 
\[
L_{+}v=(-\Delta+I)(a\cdot\nabla Q).
\]
Taking the inner product with $a\cdot\nabla Q$ yields
\[
0=\langle (-\Delta+I)(a\cdot\nabla Q), a\cdot\nabla Q\rangle,
\]
which immediately implies $a=0$, whence $L_{+}v=0$. But since $v\in {\mathcal{H}}$, the result in \cite{Lenzman} implies $v=0$. Using now the Fredholm alternative for $I+K^{*}$, we see that this operator is invertible. Thus $I+K$ is also invertible.

\medskip

In order to obtain exponential decay for $u$, we use the well-known argument of Agmon \cite{Ag}.
For $f\in{\mathcal{H}}$ and $0<c<1$ assume a bound
$
|f(x)|\lesssim e^{-c \|x\|}.
$
Write the equation above with $u\in \calH$ as 
\[
u=(-\Delta+1)^{-1}f-(-\Delta+I)^{-1}\circ K_{1}u
\]
Note that standard elliptic bootstrap methods imply $u\in H^{2}(\R^{3})$ and in particular $u\in L^{\infty}(\R^{3})$.
If we use the explicit kernel representation for $(\-\Delta +I)^{-1}$, we get the equation 
\[
u(x)=\int_{\R^{3}}\frac{e^{-\|x-y\|}}{4\pi\|x-y\|}f(y)dy -\int_{\R^{3}}\frac{e^{-\|x-y\|}}{4\pi\|x-y\|}[\phi_{Q^{2}}u+2\phi_{(Qu)}Q](y)dy,
\]
Now choose $0<\delta\leq \frac 12 c$ such that $Q(x)\lesssim e^{-\delta \|x\|}$, and pick a large $M>>1$.  We now estimate 
$\min\{M, e^{\delta \|x\|}\}u(x)$. First, we note that 
\[
 e^{\delta \|x\|} \left|\int_{\R^{3}}\frac{e^{-\|x-y\|}}{4\pi\|x-y\|}f(y)dy\right|\leq C.
\]
 Next, choose $L$ large enough such that 
\[
(\max_{|y|\geq L}|\phi_{Q^{2}}(y)|) ||\frac{e^{-(1-\delta)\|x\|}}{4\pi\|x\|}||_{L^{1}(\R^{3})}<\frac{1}{2}.
\]
Then we get 
\bea\nonumber
&&|\min\{M, e^{\delta \|x\|}\}\int_{\R^{3}}\frac{e^{-\|x-y\|}}{4\pi\|x-y\|}[\phi_{Q^{2}}u+2\phi_{Qu}Q](y)dy|\leq \frac{1}{2}\sup_{y\in\R^{3}}[\min\{M, e^{\delta \|y\|}\}|u(y)|]\\&&\nonumber\quad + 
 e^{\delta \|x\|} \left[|\int_{|y|\leq L}\frac{e^{-\|x-y\|}}{4\pi\|x-y\|} \phi_{Q^{2}}udy|
 +|\int_{\R^{3}}\frac{e^{-\|x-y\|}}{4\pi\|x-y\|}2\phi_{(Qu)}Q(y)dy|\right]
\\ &&
\nonumber \quad \leq \frac{1}{2}\sup_{y\in\R^{3}}[\min\{M, e^{\delta \|y\|}\}|u(y)|]
+C,
\eea
where $C$ is   independent of $M$. Summarizing, we have shown that 
\[
\min\{M, e^{\delta \|x\|}\}u(x)\leq C_0+\frac{1}{2}\sup_{y\in\R^{3}}[\min\{M, e^{\delta \|y\|}\}|u(y)|]
\]
with $C_0$ independent of $M$. Since $\sup_{y\in\R^{3}}[\min\{M, e^{\delta \|y\|}\}|u(y)|]<\infty$, taking the supremum over $x\in \R^{3}$ yields 
\[
\sup_{x\in\R^{3}}[\min\{M, e^{\delta \|x\|}\}|u(x)|]\leq 2C_{0},
\]
 with $C_{0}$ independent of $M$. Letting $M\rightarrow\infty$, we obtain the exponential decay of $u(x)$. 
The corresponding statement for the derivatives follows analogously, provided the derivatives of $f$ decay exponentially. 

The statement concerning the admissibility is obtained by linearity, and (ii)  is proved like (i). Parity preservation is also is standard fact.
This concludes the proof of Lemma \ref{lemmainversion}.

\medskip

Now, we proceed to the computation of the first terms of the approximate solution. This will provide an initialization of the induction argument of the proof of Proposition \ref{propinduction}. Moreover, 
the exact expression of the first two terms is required in the rest of the paper.

\bee
&&- \phi_{|v_{j+1}|^2}^{2,app} =  \frac{\lambda_j^2}{4\pi\lambda_{j+1}\|\aa\|}\int |v_{j+1}(\xi)|^2d\xi\\
& & \qquad - (-1)^{j}\frac{\lambda_j^2}{4\pi\lambda_{j+1}\|\aa \|^2}\left[\lambda_{j+1} \int \frac{(\aa\cdot\xi)}{\|\aa\|}|v_{j+1}(\xi)|^2\,  d\xi -\lambda_j \frac {(y_j\cdot\aa)}{\|\aa\|}\int|v_{j+1}(\xi)|^2d\xi\right].
\eee

\noindent{\bf Order $\frac{1}{\|\aa\|}$}:  We choose $m^{(1)}_j=b_j^{(1)}=0$ and $V_j^{(1)}=Q+T_j^{(1)}$ such that: 
$$-L_+T_j^{(1)}+ \frac{\lambda_j^2\|Q\|_{L^2}^2}{4\pi\lambda_{j+1}\|\aa\|}Q=0.$$ Note that $(Q,\nabla Q)=0$ and thus the existence of an admissible, real-valued and even function  $T_j^{(1)}$ of degree $-1$ follows from Lemma \ref{lemmainversion}.

\medskip

\noindent{\bf Order $\frac{1}{\|\aa\|^2}$}: We set $V_j^{(2)}=Q+T_j^{(1)}+T_j^{(2)}$ where $T^{(2)}_{j}$ is complex-valued and matches the $O(\frac{1}{\|\aa\|^2})$ terms. First observe from the radial symmetry of $Q$ that $\int Q^2(\xi)(\xi\cdot \aa)d\xi=0$ and thus the order 2 term in the dipole expansion is exactly $$\phi_{\|Q\|^2}^{2,2,app}=(-1)^j\
\left(\frac{\lambda_j^3\|Q\|_{L^2}^2}{4\pi\lambda_{j+1}\|\aa\|^3}\right)(\aa\cdot y_j).$$ Moreover, $T^{(1)}_j$ generates a contribution of degree $(-2)$ to the imaginary part of $\tilde{Eq}_j^{(1)}$ through the term $\frac{\partial T_j^{(1)}}{\partial \aa}\cdot2\beta$. We thus look for $m_j^{(2)}$, $b_j^{(2)}$ such that we can solve the following system:
\be
\label{eqim2}
-L_-(Im(T^{(2)}_j))+\lambda_j^2 \frac{\partial T_j^{(1)}}{\partial \aa}\cdot2\beta+\lambda_j m_j^{(2)}\Lambda Q=0,
\ee 
\bea
\label{eqre2}
\nonumber  &&   -L_+(Re(T_j^{(2)}))-\phi_{2QT_j^{(1)}}T_j^{(1)}-\phi_{(T_j^{(1)})^2}Q
\\ \nonumber&&-\frac{\lambda_j^2}{4\pi\lambda_{j+1}\|\aa\|}\left[\|Q\|_{L^2}^2 T^{(1)}_{j+1}+2\left(\int QT_{j+1}^{(1)}\right)Q\right]
- \lambda_j^3 (b_j^{(2)}\cdot y_j) Q\\
&& +(-1)^j
\left(\frac{\lambda_j^3\|Q\|_{L^2}^2}{4\pi\lambda_{j+1}\|\aa\|^3}\right)(\aa\cdot y_j) Q=0.
\eea
 From Lemma \ref{lemmainversion}, we may solve (\ref{eqim2}) with $Im(T_j^{(2)})$ admissible of order $(-2)$ provided we adjust the scaling parameter $m_j^{(2)}$ such that 
 \be
 \label{estmjdeux}
 -\lambda_j m_j^{(2)}(\Lambda Q,Q)=-\lambda_j \frac{\|Q\|_{L^2}^2}{2}m_j^{(2)}=\lambda_j^2 \left(\frac{\partial T_j^{(1)}}{\partial \aa}\cdot2\beta,Q\right) \ee so that  $m_j^{(2)}\in {\cal{S}}_2.$
 
 Similarily, we compute $b_j^{(2)}$ in order to have that the inner product of (\ref{eqre2}) with $\nabla Q$ is zero, and this yields:
 \be
\label{modulationorder2}
b^{(2)}_j=(-1)^{j+1}\frac{\|Q\|_{L^2}^2}{4\pi\lambda_{j+1}}\frac{\aa}{\|\aa\|^3}\in {\cal{S}}_2.
\ee
Note $b^{(2)}_2-b^{(2)}_1$ is the first term in the formal expansion of $\dot\beta$.  
Thus, (\ref{modulationorder2}), written for $j=1,2$ corresponds to the modulation equations on the relative position $\aa=\aa_2-\aa_1$, $\beta=\beta_2-\beta_1$: $$\aad=2\beta, \ \ \bd=-\frac{\|Q\|_{L^2}^2}{4\pi}\left(\frac{1}{\lambda_1}+\frac{1}{\lambda_2}\right)\frac{\aa}{\|\aa\|^3}$$ that is the two body problem interaction in Newtonian gravity with masses $\frac{1}{\lambda_1}, \frac{1}{\lambda_2}$.

\begin{remark}\label{rk1}
From (\ref{eqre2}) and (\ref{modulationorder2}), we infer that $Re(T_j^{(2)})$ is even, since the term in $b_j^{(2)}$ completely eliminates the only odd term in the equation of $Re(T_j^{(2)})$.
\end{remark}

\begin{remark}
Let us observe that there is no uniqueness in the choice of $T_j^{(1),(2)}$ as one can always add non trivial elements of the non trivial kernels of $L_+,L_-$. However, this would in turn modify the modulation equations and then one can explicitely check that in fact, the same {\it full} function $v_j^{(N)}(t,x)$ would be produced in the end. Hence an asymptotic expansion of this form to high order  is indeed unique.
\end{remark}


\subsection{The main induction argument}


We now proceed through the proof of Proposition \ref{propinduction} by induction on $N$.

\medskip

{\bf Proof of Proposition \ref{propinduction}}

{\bf Step 1.} Expansion of the order $N+1$ error.

We argue by induction. The case $N=2$ has been completely described in the previous subsection. We assume $N$ and prove $N+1$. Let $V_j^{(N+1)}=V_j^{(N)}+T_j^{(N+1)}$, $M_j^{(N+1)}=M_j^{(N)}+m^{(N+1)}_j$, $ B_j^{(N+1)}=B_j^{(N)}+b^{(N+1)}_j$. Then an explicit computation gives:
\bee \nonumber\tilde{Eq}_j^{(N+1)} & = & -\left[L_+(Re(T_j^{(N+1)}))+\lambda_j^3 (b_j^{(N+1)})\cdot y_j Q\right]\\
\nonumber & - &i\left[L_-(Im(T_j^{(N+1)}))+\lambda_j (m_j^{(N+1)})\Lambda Q\right]+ \tilde{Eq}_j^{(N)}-\phi^{N+1,N+1,app}_{|Q|^2}Q\\
& + & Err_j^{(1)}+iErr_j^{(2)}+Err_j^{(3)}+Err_j^{(4)}
\eee
$Err_j^{(1)}$ encodes linear type errors on the geometrical paramaters: 
\bea
\label{errorun}
\nonumber Err_j^{(1)} & = & -i\lambda_j \left[M_j^{(N)}\Lambda T_j^{(N+1)}+(m_j^{(N+1)})\Lambda \sum_{k=1}^{N+1} T_j^{(k)}\right]\\
& - & \lambda_j^3\left[B_j^{(N)}\cdot y_j T_j^{(N+1)}+(b_j^{(N+1)})\cdot y_j\sum_{k=1}^{N+1} T_j^{(k)}\right].
\eea
$Err_j^{(2)}$ is the error due to the formal time derivation: 
\bea
\label{errortwo}
\nonumber Err_j^{(2)} & & = \lambda_j^2 \sum_{k=1}^2\left[m_k^{(N+1)}\frac{\partial V_j^{(N+1)}}{\partial\lambda_k}+M_j^{(N)}\frac{\partial T_j^{(N+1)}}{\partial \lambda_j}\right]+\frac{\partial T_j^{(N+1)}}{\partial \aa}\cdot 2\beta\\
& - &  \left[\frac{\partial V_j^{(N+1)}}{\partial x}\cdot (b^{(N+1)}_{2}-b_1^{(N+1)}) +\frac{\partial T_j^{(N+1)}}{\partial \beta}\cdot (B_{2}^{(N)}-B_1^{(N)})\right].
\eea
$Err_j^{(3)}$ encodes the nonlinear type errors in the expansion of the gravitational field created  by $v_j$: 
\bea
\label{errthree}
\nonumber Err_j^{(3)}&=&\phi_{|V_j^{(N)}|^2-Q^2}T_j^{(N+1)}
+2\phi_{Re(V_j^{(N)},\overline{T_j^{(N+1)}})}(V_j^{(N+1)}-Q)\\
&+&2\phi_{Re((V_j^{(N)}-Q),\overline{T_j^{(N+1)}})}Q
+ \phi_{|T_j^{(N+1)}|^2} V_j^{(N+1)}.
\eea 
$Err_i^{(4)}$ encodes the nonlinear type errors generated by the gravitational field $\phi_i^{N,app}$: 
\bea
\label{errorfour}
\nonumber Err_j^{(4)} & = & \left(\phi^{N,app}_{|T^{(N+1)}_{j+1}|^2}+\phi^{N,app}_{2Re(V_{j+1}^{(N)},\overline{T_j^{(N+1)}})}\right)V_j^{(N+1)}\\
& + & \phi^{N+1,N+1,app}_{|Q|^2}(V_j^{N+1}-Q)+\phi^{N+1,N+1,app}_{|V_j^{(N+1)}|^2-Q^2}V_j^{(N+1)}.
\eea
From the induction hypothesis, we extract the error term of order $N+1$ of $\tilde{Eq}^{(N,n)}_j$:
$$\tilde{Eq}^{(N)}_j=\sum_{n=N+1}^{k_2(N)}\tilde{Eq}^{(N,n)}_j$$ with $\tilde{Eq}^{(N,n)}_j$ admissible of degree $n$. We then rewrite:
\bea
 \label{deferror}
&&\nonumber  \tilde{Eq}_j^{(N+1)} =  -\Big[L_+(Re(T_j^{(N+1)}))+\lambda_j^3 (b_j^{(N+1)})\cdot y_j Q\Big]\\ \nonumber &&\quad -i\Big[L_-(Im(Re(T_j^{(N+1)}))
+\lambda_j (m_j^{(N+1)})\Lambda Q\Big]
 + \tilde{Eq}^{(N,N+1)}_j-\phi^{N+1,N+1,app}_{|Q|^2}Q\\
& & \quad +Err_j^{(1)}+iErr_j^{(2)}+Err_j^{(3)}+Err_j^{(4)}+\sum_{n=N+2}^{k_2(N)}\tilde{Eq}^{(N,n)}_j
\eea

{\bf Step 2.} Computation of the order $N+1$ corrections.

Let us now compute the order $N+1$ corrections. From $(\nabla Q,yQ)=-\|Q\|_{L^2}^2$ we may chose $b_j^{(N+1)}$ such that 
$$
\ \ \left(-\lambda_j^3 b_j^{(N+1)}\cdot y_jQ+Re(\tilde{Eq}^{(N,N+1)}_j-\phi^{N+1,N+1,app}_{|Q|^2}Q),\nabla  Q\right)=0
$$From Lemma \ref{lemmainversion}, we now can solve for a $Re(T_j^{(N+1)})$ satisfying:
$$-\left[L_+(Re(T_j^{(N+1)}))+(b_j^{(N+1)}\cdot y_j ) Q\right]+Re\left[\tilde{Eq}^{(N,N+1)}_j-\phi^{N+1,N+1,app}_{|Q|^2}Q\right]=0.
$$ Similarily, from the subcritical relation $(\Lambda Q,Q)=\frac{1}{2}\|Q\|_{L^2}^2$, we chose $m_j^{(N+1)}$ satisfying 
$$
\left(-\lambda_j m_j^{(N+1)}\Lambda Q+Im(\tilde{Eq}^{(N,N+1)}_j-\phi^{N+1,N+1,app}_{|Q|^2}Q),Q\right)=0
$$and thus we can solve for: 
$$
-\left[L_-(Im(T_j^{(N+1)}))+(m_j^{(N+1)})\Lambda Q\right]+Im\left[\tilde{Eq}^{(N,N+1,m)}_j-\phi^{N+1,N+1,app}_{|Q|^2}Q\right]=0.
$$
Observe now that $\phi^{N+1,N+1,app}_{|Q|^2}Q$ is admissible of degree $N+1$ and thus from the induction hypothesis and Lemma \ref{lemmainversion}, we just count the homogenity in $\aa$ and find that $T_j^{(N+1)}$ is admissible of degree $N+1$ and $m_j^{(N+1)}\in {\cal{S}}_{N+1}$, $b_j^{(N+1)}\in {\cal{S}}_{N+1}$.

\medskip

{\bf Step 3.} Estimate of the error.

We now estimate the error terms in the right side of (\ref{deferror}). Note that the last term in the RHS of (\ref{deferror}) is by constrution admissible of order $\geq N+2$. It thus remains to consider the $(Err_j^{(k)})_{1\leq k\leq 4}$.
\begin{itemize}
\item $Err_j^{(1)}$: Consider (\ref{errorun}). From the construction of $v_j^{(N+1)}, M_j^{(N+1)}$ and  $B_j^{(N+1)}$, and from $deg(V_j^{(N+1)}-Q)\geq 1,$ $deg(M_j^{(N)})\geq 2, deg(B_j^{(N)})\geq 2$ and $deg(T_j^{(N+1)})$, $deg(m_j^{(N+1)}),$ $deg(b_j^{(N+1)})\geq N+1$,  $Err_j^{(1)}$ is admissible of degree $\geq N+2$.
\item $Err_j^{(2)}$: Consider (\ref{errortwo}). We observe that given $\psi$ admissible of degree $n$, $\frac{\partial\psi}{\partial \lambda}, \frac{\partial\psi}{\partial \beta}$ are admissible of degree $n$ while $\frac{\partial\psi}{\partial x}$ is admissible of degree $n+1$. We then use the fact that $deg(\frac{\partial V_j^{(N+1)}}{\partial\lambda_k})\geq 3$ and the key degeneracy  $deg(B_j^{(N)})\geq 2$ to conclude that $Err_j^{(2)}$ is admissible of degree $\geq N+2$.
\item $Err_j^{(3)}$: Consider (\ref{errthree}). Then $deg(|V_j^{(N)}|^2-Q^2)\ge 1$ and $deg(V_j^{(N+1)}-Q)\geq 1$ together with Lemma \ref{stabadmissible} imply that $Err_j^{(3)}$ is admissible of degree $\geq N+2$.
\item $Err_j^{(4)}$: Consider (\ref{errorfour}). The key is to recall from (iii) of Lemma \ref{stabadmissible} that $\phi^{N,app}_{\psi}$ always adds at least one to the degree of $\psi$. Using again $deg(V_j^{(N+1)}-Q)\geq 1$, we conclude that $deg(Err_j^{(3)})\geq N+2$.
\end{itemize}

This concludes the induction step and the proof of Proposition \ref{propinduction}.


\subsection{Definition of the approximate solution}


Now we explicitely describe how to extract a high order approximate solution to (\ref{hartree}) from the profiles defined in Proposition \ref{propinduction}.

\begin{Prop}[Definition and accuracy of the approximation]
\label{constructionansatz}
Given $N\geq 1$ and parameters $P=(\aa_1,\aa_2,\beta_1,\beta_2, \lambda_1,\lambda_2)\in (\R^3)^4\times (\R^*_+)^2\times(\R)^2$, let
$$(V_j^{(N)})_P(y_j)=V_j^{(N)}(\aa_2-\aa_1,\beta_2-\beta_1,\lambda_1,\lambda_2, y_j),\ \ j=1,2$$ be the object constructed from Proposition \ref{propinduction}.

Let $t\to P(t)=(\aa_1(t),\aa_2(t),\beta_1(t),\beta_2(t), \lambda_1(t),\lambda_2(t))$ and $\gamma(t)=(\gamma_1(t),\gamma_2(t))$ be  ${\cal{C}}^1$ functions on a time interval $I=[t_0,t_1]$, $t_1\leq +\infty$. Let $\aa(t)=\aa_2(t)-\aa_1(t)$ and $\beta(t)=\beta_2(t)-\beta_1(t)$ and assume that the following controls hold: $\forall t\in I=[t_0,t_1]$, 
\be
\label{controlsxb}
1\leq\frac{\|\aa(t_0)\|}{2}\leq \|\aa (t)\|, \ \ \|\beta(t)\|\leq 2\|\beta(t_0)\|,
\ee
\be
\label{controlscaling}
0<\frac{\lambda_j(t_0)}{2}\leq \lambda_j(t)\leq 2\lambda_j(t_0), \ \ j=1,2.
\ee
Define the approximate two soliton solution by: 
\bea
\label{defR}
R^{(N)}(t,x)&  =&    R^{(N)}_{P(t),\gamma(t)}(x)= R^{(N)}_1(t,x)+R^{(N)}_2(t,x)\\
\nonumber & = & \sum_{j=1}^2 \frac{1}{\lambda_j^2(t)}(V_j^{(N)})_{P(t)}\left(\frac{x-\aa_j(t)}{\lambda_j(t)}\right)e^{-i\gamma_j(t)+i\beta_j(t)\cdot x}.
\eea
Then $R^{(N)}$ satisfies
\be
\label{eqR}
 i\partial_tR^{(N)} + \Delta R^{(N)}-\phi_{|R^{(N)}|^2}R^{(N)}=-\Psi^{(N)}+\sum_{j=1}^2\frac{1}{\lambda_j^4}\left[S_j^{(N)}\right]\left(\frac{x-\aa_j(t)}{\lambda_j(t)}\right)e^{-i\gamma_j(t)+i\beta_j(t)\cdot x}
\ee
where $S_j^{(N)}$ given by (\ref{defsnj}) encodes the dynamical system of the geometrical parameters 
and $\Psi^{(N)}$ collects the error terms and satisfies:
\be
\label{sizeerror}
\forall x \in \R^N, \ \ |\Psi^{(N)}(t,x)|\leq \frac{1}{\|\aa(t)\|^{N+1}}e^{-C_{N,\lambda_1,\lambda_2,\beta_1,\beta_2}(\|x-\aa_1(t)\|+\|x-\aa_2(t)\|)}.
\ee
\end{Prop}

 \textbf{Proof.} According to the computations leading to the derivation of the approximate equation (\ref{approximationone}), there holds the explicit expression of $\Psi^{(N)}$:
\bea
\label{expressionerror}
& &  \Psi^{(N)} = 2\phi_{Re(R_1^{(N)}\overline{R_2^{(N)}})}(R_1^{(N)}+R_2^{(N)}) \\
\nonumber &+&  \sum_{j=1}^2 \frac{1}{\lambda_j^4}\left[\tilde{Eq}^{(N)}_j\right]\left(\frac{x-\aa_j(t)}{\lambda_j(t)}\right)e^{-i\gamma_j(t)+i\beta_j(t)\cdot x}+ e^{-i\gamma_j(t)+i\beta_j(t)\cdot x}
\\
\nonumber &  &   \sum_{j=1}^2\frac{1}{\lambda_j^4}\left\{\left[\phi^{N,app}_{|V^{(N)}_{j+1}|^2}-\left(\frac{\lambda_j}{\lambda_{j+1}}\right)^2\phi_{|V^{(N)}_{j+1}|^2}\left(\frac{\lambda_j}{\lambda_{j+1}}y_{j}{+}(-1)^j\frac{\aa}{\lambda_{j+1}}\right)\right]V_j^{(N)}\right\}\left(\frac{x{-}\aa_j(t)}{\lambda_j(t)}\right)\eea
with $\tilde{Eq}^{(N)}_j$ given by (\ref{defeqtilden}).

The first term in the RHS of (\ref{expressionerror}) is estimated from the exponential localization in space of the profile $R^{(N)}_j$:
\bee
\left|2\phi_{Re(R_1^{(N)}\overline{R_2^{(N)}})}(R_1^{(N)}+R_2^{(N)})\right| & \leq & C\|\phi_{Re(R_1^{(N)}\overline{R_2^{(N)}})}\|_{L^{\infty}}(|R_1^{(N)}|+|R_2^{(N)}|)\\& \leq & \|R_1^{(N)} R^{(N)}_2\|_{L^{\frac{3}{2}}}e^{-C_{N,\lambda_1,\lambda_2,\beta_1,\beta_2}(\|x-\aa_1(t)\|+\|x-\aa_2(t)\|)}\\
& \leq & e^{-C_{N,\lambda_1,\lambda_2,\beta_1,\beta_2}(\|\aa(t)\|)+\|x-\aa_1(t)\|+\|x-\aa_2(t)\|}
\eee 
where we used the standard Hardy-Littlewood-Sobolev inequalities.

According to Proposition \ref{propinduction}, the term $\tilde{Eq}^{(N)}_j$ in the RHS of (\ref{expressionerror}) is admissible of degree $\geq(N+1)$ and estimate (\ref{sizeerror}) follows for this term.

It remains to estimate the third term in the RHS of (\ref{expressionerror}) that is the error induced by the dipole expansion of the gravitational field created by $V^{(N)}_{j+1}.$ We claim that:
\bea
\label{errordipole}
\nonumber & & \left|\phi^{N,app}_{|V^{(N)}_{j+1}|^2}(y_j)-\left(\frac{\lambda_j}{\lambda_{j+1}}\right)^2\phi_{|V^{(N)}_{j+1}|^2}\left(\frac{\lambda_j}{\lambda_{j+1}}y_{j}+(-1)^j\frac{\aa}{\lambda_{j+1}}\right)\right|\\
& \leq & C_{N,\lambda_1,\lambda_2,\beta_1,\beta_2}(1+\|y_j\|)^{N+1}\left[\frac{{\bf 1}_{\lambda_j\|y_j\|\leq \frac{\|\aa\|}{10}}}{\|x\|^{N+1}}+{\bf 1}_{\lambda_j\|y_j\|\geq \frac{\|\aa\|}{10}}\right].
\eea
The exponential localization of $V_j^{(N)}$ in $y_j$ then implies: $$\left|\phi^{N,app}_{|V^{(N)}_{j+1}|^2}-\left(\frac{\lambda_j}{\lambda_{j+1}}\right)^2\phi_{|V^{(N)}_{j+1}|^2}\left(\frac{\lambda_j}{\lambda_{j+1}}y_{j}+(-1)^j\frac{\alpha}{\lambda_{j+1}}\right)\right|V_j^{(N)}\leq \frac{e^{-C_{N,\lambda_1,\lambda_2,\beta_1,\beta_2}\|y_j\|}}{\|\alpha\|^{N+1}}$$ which concludes the proof of (\ref{sizeerror}).

\medskip

Proof of (\ref{errordipole}): We first observe from exponential localization of $V_j^{(N+1)}$ and the definition (\ref {approxfield}) of the dipole expansion that: $\forall y_j\in \R^N$, 
$$ \left|\phi^{N,app}_{|V^{(N)}_{j+1}|^2}(y_j)\right|+\left|\left(\frac{\lambda_j}{\lambda_{j+1}}\right)^2\phi_{|V^{(N)}_{j+1}|^2}\left(\frac{\lambda_j}{\lambda_{j+1}}y_{j}+(-1)^j\frac{\alpha}{\lambda_{j+1}}\right)\right|\leq  C_{N,\lambda_1,\lambda_2,\beta_1,\beta_2}(1+\|y_j\|)^{N+1}.$$ We may thus restrict our attention to the region $\lambda_j\|y_j\|\leq\frac{\|\alpha\|}{10}.$ We then split:
\bee
& &   \left(\frac{\lambda_j}{\lambda_{j+1}}\right)^2\phi_{|V^{(N)}_{j+1}|^2}\left(\frac{\lambda_j}{\lambda_{j+1}}y_{j}+(-1)^i\frac{\alpha}{\lambda_{j+1}}\right)\\ & =& -\frac{\lambda_j^2}{4\pi\lambda_{j+1}}\int_{\R^3}\frac{|V^{(N)}_{j+1}(\xi)|^2}{\|\alpha-(-1)^j(\lambda_{j+1}\xi-\lambda_j y_j)\|}d\xi\\
 & = & -\frac{\lambda_j^2}{4\pi\lambda_{j+1}}\left[\int_{\lambda_{j+1}\|\xi\|\leq \frac{\|\aa\|}{10}}\frac{|V^{(N)}_{j+1}(\xi)|^2}{\|\aa-(-1)^j(\lambda_{j+1}\xi-\lambda_jy_j)\|}d\xi\right .\\
 & + & \left . \int_{\lambda_{j+1}\|\xi\|\geq \frac{\|\aa\|}{10}}\frac{|V^{(N)}_{j+1}(\xi)|^2}{\|\aa-(-1)^j(\lambda_{j+1}\xi-\lambda_jy_j)\|}d\xi\right]
 \eee
 From the exponential localization in space of $V^{(N)}_{j+1}$ and Hardy-Littlewood-Sobolev: 
 \bee
 \left\|\int_{|\xi\|\geq \frac{\|\aa\|}{10}}\frac{|V^{(N)}_{j+1}(\xi)|^2}{\|\aa-(-1)^j(\lambda_{j+1}\xi-\lambda_jy_j)\|}d\xi\right\|_{L^{\infty}} & \leq  & C_{\lambda_1,\lambda_2}\left\|V^{(N)}_{j+1}{\bf 1}_{\|\xi\|\geq \frac{\|\aa\|}{10}}\right\|_{L^{\frac{8}{3}}}^2\\
 & \leq & e^{-C_{N,\lambda_1,\lambda_2,\beta_1,\beta_2}\|\aa\|}.
 \eee
Now for $\lambda_j\|y_j\|\leq\frac{\|\aa\|}{10}$ and $\lambda_{j+1}\|\xi\|\leq \frac{\|\aa\|}{10}$, we use (\ref{expansionofthefield}) that is the uniform estimate:
$$\forall \|\aa\|\geq 1, \ \ \forall \|\zeta\|\leq \frac{\|\aa\|}{2}, \ \ \left|\frac{1}{\|\aa-\zeta\|}-\sum_{k=1}^NF_k(\aa,\zeta)\right|\leq C_N\frac{(1+\|\zeta\|)^{N+1}}{\|\aa\|^{N+1}}$$ to derive: 
\bee
& & \left| -\frac{\lambda_j^2}{4\pi\lambda_{i+1}}\left[\int_{\lambda_{j+1}\|\xi\|\leq \frac{\|\aa\|}{10}}\frac{|V^{(N)}_{j+1}(\xi)|^2}{\|\aa-(-1)^j(\lambda_{j+1}\xi-\lambda_jy_j)\|}d\xi\right]-\phi^{N,app}_{|V^{(N)}_{j+1}|^2}(y_j)\right|\\
& \leq & C_{N,\lambda_1,\lambda_2,\beta_1,\beta_2}\frac{(1+\|y_j\|)^{N+1}}{\|\aa\|^{N+1}}\int_{\R^3}|V^{(N)}_{j+1}(\xi)|^2(1+\|\xi\|)^{N+1}d\xi
\eee
and (\ref{errordipole}) from the exponential localization of $V_{j+1}^{(N)}$. This concludes the proof of (\ref{errordipole}).\\
This concludes the proof of Proposition \ref{constructionansatz}.


\section{Reduction of the problem}

In this section, we reduce the proof of Theorem \ref{thmhyperbolic} to a bootstrap argument i.e. Lemma~\ref{boot}. The reduction proceeds into two  steps: a compactness argument assuming a uniform estimate (section \ref{sec:3.1}) and a proof of the uniform estimate assuming that the bootstrap holds (section \ref{sec:3.2}).

This approach follows  the strategy developped in \cite{Mcrit}, \cite{Mkdv} and allows us to use an energy type method for the construction.

\medskip

Let $(\lambda^{\infty}_1,\lambda_2^{\infty})\in(\R^*_+)^2$ and a solution $$P^{\infty}(t)=(\aa_1^\infty(t),\aa_2^\infty(t),\beta_1^\infty(t),\beta_2^\infty(t),\lambda_1^\infty,\lambda_2^\infty)$$ of the Newtonian two body problem:
\be
\label{defpinfty}
  \left   \{ \begin{array}{llll}
         \aad^{\infty}_1=2\beta^{\infty}_1, \ \ \aad^{\infty}_2=2\beta^{\infty}_2, \ \ \aa^{\infty}=\aa_2^{\infty}-\aa_1^{\infty},\\
         \bd^{\infty}_1=\frac{\|Q\|_{L^2}^2}{4\pi\lambda_2^\infty}\frac{\aa^{\infty}}{\|\aa^{\infty}\|^3}, \ \ \bd^{\infty}_2=-\frac{\|Q\|_{L^2}^2}{4\pi\lambda_1^\infty}\frac{\aa^{\infty}}{\|\aa^{\infty}\|^3},         \end{array}
\right .
\ee
satisfying 
\be
\label{cnscnoson}
\lambda^{\infty}_2 \aa^{\infty}_1(t)+\lambda^{\infty}_1 \aa^{\infty}_2(t)=0,\quad
(\aa_{j})_{(3)}^\infty(t)=0,
\ee
and either of hyperbolic type with any masses 
$$E_0=\frac 14 \|\aad^{\infty}\|^2-\frac{\|Q\|_{L^2}^2}{4\pi}\left(\frac{1}{\lambda^{\infty}_1}+\frac{1}{\lambda^{\infty}_2}\right)\frac{1}{\|\aa^{\infty}\|}> 0
$$ 
or of parabolic type with the restriction of equal masses:
$$E_0=0\ \ \mbox{and}\ \ \lambda_1^{\infty}=\lambda_2^{\infty}=\lambda^{\infty}.
$$

Recall that  the two particles move away from each other with a speed depending on the regime
\be
\label{estxhyp}
\mbox{Hyperbolic regime:}\qquad
\|\aa^{\infty}(t)\|\mathop{\sim}_{+\infty} C_{P^{\infty}}t, \ \ \|\beta^{\infty}(t)\|\mathop{\sim}_{+\infty} C>0,
\ee
\be
\label{estxparab}
\mbox{Parabolic regime:}\qquad
\|\aa^{\infty}(t)\|\mathop{\sim}_{+\infty} C_{P^{\infty}}t^{\frac{2}{3}}, \ \ \|\beta^{\infty}(t)\|\mathop{\sim}_{+\infty}\frac{C}{t^{\frac{1}{3}}}.
\ee
By rotation invariance, in addition to the assumption that the $P^{\infty}$ trajectory takes place in the plane $\{x_{(3)}=0\}$, we may also require without loss of generality that  in this plane, the line $\{x_{(1)}=0\}$ strictly separates the two particles, in the following sense
\be
\label{behaviorx}
\liminf_{t\to +\infty} \frac {\aa_1^\infty(t)} {t^{\mu}} \cdot \vec{e}_1 >0, 
\ \ \liminf_{t\to +\infty} \frac {\aa_2^\infty(t)} {t^{\mu}} \cdot (-\vec{e}_1) >0,
\ee
where $\mu=1$ or $\frac 23$, depending on the regime.


\subsection{A compactness argument}\label{sec:3.1}

 
 Let us first show that the system of modulation equations inherited from Proposition \ref{propinduction} and (\ref{defsnj}) may be solved from infinity with asymptotic behavior given by $P^{\infty}$.

\begin{Lemma}[Integration from infinity of the refined modulation equations]
\label{defdynamic}\quad \\ Let $M_j^{(N)}, B_j^{(N)}$ be given by (\ref{computmb}) in Proposition \ref{propinduction}. There exists  a solution $$P^{(N)}=(\aa_1^{(N)}(t),\aa_2^{(N)}(t),\beta_1^{(N)}(t),\beta_2^{(N)}(t),\lambda_1^{(N)}(t),\lambda_2^{(N)}(t))$$ 
of the system
\begin{equation}\label{bofbof} \left   \{ \begin{array}{lllll}
	\ld^{(N)}_1=M_1^{(N)}(\aa^{(N)},\beta^{(N)},\lambda_1^{(N)},\lambda_2^{(N)}) \ \ \ld^{(N)}_2=M_2^{(N)}(\aa^{(N)},\beta^{(N)},\lambda_1^{(N)},\lambda_2^{(N)}) ,\\
         \aad^{(N)}_1=2\beta^{(N)}_1, \ \ \aad^{(N)}_2=2\beta^{(N)}_2, \ \ \aa^{(N)}=\aa_2^{(N)}-\aa_1^{(N)},\ \ \beta^{(N)}=\beta_2^{(N)}-\beta_1^{(N)},\\
         \bd^{(N)}_1=B_1^{(N)}(\aa^{(N)},\beta^{(N)},\lambda_1^{(N)},\lambda_2^{(N)}), \ \ \bd_2^{(N)}=B_2^{(N)}(\aa^{(N)},\beta^{(N)},\lambda_1^{(N)},\lambda_2^{(N)}), 
         \end{array}
\right .
\end{equation}
on $[T_0,+\infty)$ for $T_0$ large enough, satisfying in addition:
\begin{enumerate}
\item Hyperbolic regime: $$
\|P^{\infty}(t)-P^{(N)}(t)\|\leq  \frac{C_{P^{\infty}}}{t^{\frac{1}{2}}}.$$ 
\item Parabolic regime: Assume $\lambda_1^{\infty}=\lambda_2^{\infty}$. Then, $\forall j=1,2$,
$$
\sum_{k=1}^2\left |\frac{(\aa_j^{\infty})_{(k)}(t)}{(\aa^{(N)}_j)_{(k)}(t)}-1\right|+|(\aa^{(N)}_j)_{(3)}(t)|+\|\beta^{(N)}_j(t)-\beta_j^{\infty}(t)\|+|\lambda^{(N)}_j(t)-\lambda_j^{\infty}|\leq \frac{C_{P^{\infty}}}{t^{\frac{1}{3}}},
$$
$$\lambda^{(N)}_1(t)\aa^{(N)}_2(t)+\lambda^{(N)}_2(t)\aa^{(N)}_1(t)\to 0 \ \ \mbox{as} \ \ t\to +\infty.$$

\end{enumerate}
\end{Lemma}

{\bf Proof of Lemma \ref{defdynamic}.} See Appendix A.

\begin{remark} In particular, there holds using (\ref{estxhyp}), (\ref{estxparab}) as $t\to +\infty$:
\be
\label{estxbhyp}
\|\aa^{(N)}(t)\|\mathop{\sim}_{+\infty}C_{P^{\infty}}t, \ \ \|\beta^{(N)}(t)\|\mathop{\sim}_{+\infty} C_{P^{\infty}}>0\ \ \mbox{in the hyperbolic case},
 \ee
\be
\label{estxbpar}\|\aa^{(N)}(t)\|\mathop{\sim}_{+\infty} C_{P^{\infty}}t^{\frac{2}{3}}, \ \ \|\beta_1^{(N)}(t)\|+\|\beta_2^{(N)}(t)\|\leq \frac{C_{P^{\infty}}'}{t^{\frac{1}{3}}} \ \mbox{in the parabolic case}.
\ee
 \end{remark}

We claim the following uniform backwards estimate, which is the heart of the proof of Theorem \ref{thmhyperbolic}.

\begin{Prop}[Uniform estimates]
\label{propuniform}
For $N>1$ large enough the following holds.
Let $P^{(N)}(t)$ be the solution obtained from Lemma \ref{defdynamic} and let $\gamma^{(N)}(t)=(\gamma_1^{(N)}(t),\gamma_2^{(N)}(t))$
satisfy
$$
j=1,2,\ \ \dot\gamma_j^{(N)}=-\frac 1 {(\lambda_j^{(N)})^2}+\|\beta_j^{(N)}\|^2 
- \dot \beta_j^{(N)} \cdot \aa_j^{(N)}.
$$
Let 
$R^{(N)}_{P^{(N)},\gamma^{(N)}}(t)$ be defined by (\ref{defR}).

Let a sequence $T_n\to +\infty$ and $u_n(t)$ be the solution to 
\be
\label{hartreeTn}
 \ \  \left   \{ \begin{array}{ll}
         i\partial_tu_n+ \Delta u_n-\phi_{|u_n|^2}u_n=0,\ \ \phi_{|u_n|^2}=-\frac{1}{4\pi \|x\|}\star |u_n|^2,\\
        u_n(T_n)=R^{(N)}_{P^{(N)}}(T_n).
         \end{array}
\right .
\ee
Then there exists $T_0=T_0(N)$ independent of $n$ such that:
\be
\label{uniformt0}
\forall n\geq 1, \  \ \forall t\in [T_0,T_n], \ \ \left\|u_n(t)-R^{(N)}_{P^{(N)}(t)}\right\|_{H^1}\leq \frac{1}{t^{\frac{N}{8}}}.
\ee
\end{Prop}

Let us prove Theorem \ref{thmhyperbolic}  as a direct consequence of Proposition \ref{propuniform}.

\medskip

{\bf Proof of Theorem \ref{thmhyperbolic}}. This proof follows the lines of the one of   Theorem 1 in \cite{MMTnls}.
By time translation invariance, we suppose without loss of generality that $T_0=0$ in Proposition \ref{propuniform}.

By (\ref{uniformt0}), it is clear that there exists $C>0$ such that for any $n$, 
$t\in [0,T_n]$,
\be\label{unifest}
\|u_n(t)\|_{H^1}\leq C.
\ee

Now, we claim the following compactness result.

\begin{Lemma}\label{lem2}
There exists $U_0\in H^1(\R^3)$ and a subsequence $(u_{n'})$ of $(u_n)$ such that
$$
u_{n'}(0)\rightarrow U_0\hbox{ in $L^2(\R^3)$ as $n\rightarrow +\infty$.}
$$
\end{Lemma}
We  sketch the proof of Lemma \ref{lem2}. Note that local convergence in $L^2$ for a subsequence is obvious by the $H^1$ bounds. The lemma thus follows from the following uniform $L^2$ localization property
\be\label{compl2}
\forall \epsilon_0>0, \ \exists A_0=A_0(\epsilon_0),\ \mbox{such that} \ \forall n\geq 1,\
\int_{\|x\|>A_0} |u_n(0,x)|^2 dx < \epsilon.
\ee
Note finally that (\ref{compl2}) is a consequence of (\ref{uniformt0}) (here $T_0=0$) and simple computations. See the proof of Lemma 2 in \cite{MMTnls} for more details.

\medskip

Consider now the global $H^1$ solution $U(t)$ of
\be
\label{hartreeU}
 \left   \{ \begin{array}{lll}
         iU_t + \Delta U-\phi_{|U|^2}U=0,\quad 
         \phi_{|U|^2}=-\frac{1}{4\pi \|x\|}\star |U|^2, \\
         (t,x)\in \R\times \R^3, \  \ U(0)=U_0.
         \end{array}
\right .
\ee
Fix $t\geq 0$. For $n$ large enough, we have $T_n\geq t$ and by continuous dependence
of the solution of (\ref{hartree}) on the initial data in $L^2(\R^3)$
(see \cite{Cbook}), we have
$$
u_{n'}(t)\rightarrow U(t)\hbox{ in $L^2(\R^3)$ as $n\rightarrow +\infty$.}
$$
Since $u_{n'}(t)-R^{(N)}_{P^{(N)}(t)}$ is uniformly bounded in $H^1$, it follows that
$$
u_{n'}(t)-R^{(N)}_{P^{(N)}(t)}\rightharpoonup
 U(t)-R^{(N)}_{P^{(N)}(t)}\hbox{ in $H^1(\R^3)$ as $n\rightarrow +\infty$,}
$$
and thus, by (\ref{uniformt0}), we obtain
\be
\label{g10}
\ \forall t\geq 0, \ \ \|U(t)-R^{(N)}_{P^{(N)}(t)}\|_{H^1}\leq \frac{1}{t^{\frac{N}{8}}}.
\ee
Therefore, Theorem \ref{thmhyperbolic} follows from the definition of $R^{(N)}_{P^{(N)}(t)}$
in Proposition \ref{constructionansatz} and the properties of ${P^{(N)}(t)}$ described in Lemma \ref{defdynamic}.


\subsection{The bootstrap argument}\label{sec:3.2}


We now focus onto the proof of Proposition \ref{propuniform} which relies on a bootstrap argument based on energy type estimates.

\medskip

{\bf Proof of Proposition \ref{propuniform}.}
Let $N$ large enough to be chosen and let $R^{(N)}_{P^{(N)}(t),\gamma^{(N)}(t)}$ be chosen as in the statement of Proposition \ref{propuniform}.

Let $u_n(t)$ be the solution to (\ref{hartreeTn}). As long as the flow evolves close to the two solitary wave flow, we may from standard arguments (see e.g. \cite{MMnls}, Lemma 3) use the implicit function theorem to prove the existence of a unique ${\cal{C}}^1$ geometrical decomposition of the flow:
$$u_n(t,x)=R^{(N)}_{P(t),\gamma(t)}(x)+\e(t,x)$$ (see (\ref{defR})) with 
$$P(t)=(\aa_1(t),\aa_2(t),\beta_1(t),\beta_2(t), \lambda_1(t),\lambda_2(t)),\quad \gamma(t)=(\gamma_1(t),\gamma_2(t))$$
 chosen such that the following orthogonality conditions hold: for $j=1,2$, let
\be
\label{defej}
\e_j(t,y_j)=\lambda_j^2(t)\e(t,\lambda_j(t)y_j+\aa_j(t))e^{i\gamma_j(t)-i\beta_j(t)\cdot(\lambda_j(t)y_j+\aa_j(t))},
\ee
then:
\be
\label{orthconditions}
Re(\e_j(t),\overline{Q})=Re(\e_j(t), \overline{y_jQ})=Im(\e_j(t),\overline{\Lambda Q})=Im(\e_j(t),\overline{\nabla Q})=0.
\ee
Observe from (\ref{hartreeTn}) that: 
\be
\label{initialization}
P(T_n)=P^{(N)}(T_n), \ \ \gamma(T_n)=\gamma^{(N)}(T_n), \ \ \e(T_n)=0.
\ee
We claim the following bootstrap lemma which is the core of the argument:

\begin{Lemma}[Bootstrap]
\label{boot}
For $N$ large enough, there exists $T_0(N)>0$ large enough such for all $T^*\in[T_0,T_n]$, if:
\be
\label{hypflowboot}
\forall t\in[T^*,T_n],   \ \  \left   \{ \begin{array}{lll}

	 \|\e(t)\|_{H^1}\leq \frac{1}{t^{\frac{N}{4}}},\\
	 \sum_{j=1}^2\left[|\lambda_j(t)-\lambda_j^{(N)}(t)|+\|\beta_j-\beta_j^{(N)}\|\right]\leq\frac{1}{t^{1+\frac{N}{8}}},\\
	 \sum_{j=1}^2\left[|\gamma_j(t)-\gamma_j^{(N)}(t)|+\|\aa_j-\aa_j^{(N)}\|\right]\leq\frac{1}{t^{\frac{N}{8}}},
         \end{array}
\right .
\ee
then:
\be
\label{hypflowboost}
\forall t\in[T^*,T_n],   \ \  \left   \{ \begin{array}{lll}

	 \|\e(t)\|_{H^1}\leq \frac{1}{2t^{\frac{N}{4}}},\\
	 \sum_{j=1}^2\left[|\lambda_j(t)-\lambda_j^{(N)}(t)|+\|\beta_j-\beta_j^{(N)}\|\right]\leq\frac{1}{2t^{1+\frac{N}{8}}},\\
	 \sum_{j=1}^2\left[|\gamma_j(t)-\gamma_j^{(N)}(t)|+\|\aa_j-\aa_j^{(N)}\|\right]\leq\frac{1}{2t^{\frac{N}{8}}}.
         \end{array}
\right .
\ee
\end{Lemma}

From a standard continuity argument, Lemma \ref{boot} proves that (\ref{hypflowboost}) is indeed satisfied on $[T_0,T_n]$.
Thus, $\forall t\in [T_0,T_n]$, 
$$\|u_n(t)-R^{(N)}_{P^{(N)}(t),\gamma^{(N)}(t)}\|\leq \|u_n(t)-R^{(N)}_{P(t),\gamma(t)}\|+\|R^{(N)}_{P(t),\gamma(t)}-R^{(N)}_{P^{(N)},\gamma^{(N)}(t)}\|\leq \frac{C}{t^{\frac{N}{4}}}$$ which implies (\ref{uniformt0}) on $[T_0,T_n]$  and concludes the proof of Proposition \ref{propuniform}.

\medskip

The rest of the paper is devoted to the proof of  Lemma \ref{boot}.


\section{Control of the dynamics in the bootstrap regime}


This section is devoted to the proof of   Lemma \ref{boot} which relies on energy type -- and not dispersive -- estimates on $\e$.


\subsection{Control of the modulation parameters}


We focus onto the proof of the control of the modulation parameters (\ref{hypflowboost}) which is a consequence of our specific choice of orthogonality conditions (\ref{orthconditions}). 

\medskip

We follow the notation of the proof of Proposition \ref{propuniform}. We simplify the notation
$R(t)=R^{(N)}_{P(t),\gamma(t)}=R_1(t)+R_2(t)$. First observe that (\ref{estxbhyp}), (\ref{estxbpar}) and the bootstrap assumption (\ref{hypflowboot}) imply on $[T_0,T_N]$ for $T_0$ large enough:
\be
\label{estintermhyp}
\|\aa(t)\|\mathop{\sim}_{+\infty} C_{P^{\infty}}t, \ \ \|\beta(t)\|\mathop{\sim}_{+\infty} C_{P^{\infty}}>0 \ \ \mbox{in the hyperbolic case},
\ee
\be
\label{estintermpar}
\|\aa(t)\|\mathop{\sim}_{+\infty} C_{P^{\infty}}t^{\frac{2}{3}}, \ \ \|\beta_1(t)\|+\|\beta_2(t)\|\leq \frac{C'_{P^{\infty}}}{t^{\frac{1}{3}}}\ \ \mbox{in the parabolic case}.
\ee

{\bf Step 1.} Structure of the modulation equations.

Let us denote by $Mod(t)$ the errors of the modulation equations with respect to the law predicted by the dynamical system $P^{(N)}$:
\be
\label{mod}
Mod(t)=\sum_{j=1}^2 \|\aad_j-2\beta_j\|+\|\ld_j-M_j^{(N)}\|+\|\bd_j-B_j^{(N)}\|+|\gd_j+\frac{1}{\lambda_j^2}-\|\beta_j\|^2-\bd_j\cdot \aa_j|
\ee
where $M_j^{(N)}=M_j^{(N)}(P)=M_j^{(N)}(\aa, \beta,\lambda_1, \lambda_2)$,
$B_j^{(N)}=B_j^{(N)}(P)=B_j^{(N)}(\aa, \beta,\lambda_1, \lambda_2)$.

Now, we write the equation for $\e$. From (\ref{eqR}) and $\e=u_n-R$, there holds:
\bea
\label{eqe}
i\partial_t\e+\Delta \e-2\phi_{Re(\e \overline{R})}R-\phi_{|R^2|}\e & = & \Psi^{(N)}+\calN(\e)\\
\nonumber & - & \sum_{j=1}^2\frac{1}{\lambda_j^4}\left[S_j^{(N)}\right](\frac{x-\aa_j(t)}{\lambda_j(t)})e^{-i\gamma_j(t)+i\beta_j(t)\cdot x}
\eea
 with 
 \be
 \label{defne}
 \calN(\e)=\phi_{|\e|^2}R+2\phi_{Re(\e \overline{R})}\e+\phi_{|\e|^2}\e.
 \ee
Next, we compute the modulation equations using the orthogonality conditions (\ref{orthconditions}). It suffices to take the inner product of (\ref{eqe}) by four given functions exponentially localized around each center of mass $(\aa_j(t))_{j=1,2}$.
Let $\Theta_j=\Theta_j(x)$ be a given function with some decay at infinity and let
 $$\Theta(t,x)=\frac{1}{\lambda_j^2}\Theta_j\left(\frac{x-\aa_j(t)}{\lambda_j}\right)e^{-i\gamma_j+i\beta_j\cdot x},  \ \ |\Theta(t,x)|\leq e^{-C\|x-\aa_j(t)\|}.$$ 
 We compute from (\ref{eqe}) after an integration by parts:
 \bee
&& \frac{d}{dt} \left\{Im\int\e\overline{\Theta}\right\}= Re\int\e\left[\overline{i\partial_t\Theta+\Delta \Theta-2\phi_{Re(\Theta\overline{R})}R-\phi_{|R^2|}\Theta}\right]\\ 
\nonumber & & - Re(\int [\Psi^{(N)}+\calN(\e)]\overline{\Theta})+ \sum_{j=1}^2Re\left(\int\frac{1}{\lambda_j^4}\left[S_j^{(N)}\right](y_j)e^{-i\gamma_j(t)+i\beta_j(t)\cdot x}\overline{\Theta}\right).
\eee
We now use the fact that the approximate modulation equations for scaling and Galilei are both of order $O(\frac{1}{\|\aa\|^2})$ and the localization of $\Theta$ around $\aa_j$ to derive:
\bee
& &  i\partial_t\Theta+\Delta \Theta-2\phi_{Re(\Theta\overline{R})}R-\phi_{|R^2|}\Theta=\frac{1}{\lambda_j^4}\left[i\lambda_j^2\partial_t\Theta_j-L_j\Theta_j\right](y_j)e^{-i\gamma_j(t)+i\beta_j(t)\cdot x}\\
& -  & 2\phi_{Re(\Theta\overline{R_j})}R_{j+1}+O\left(\frac{1}{\|\aa\|^2}+Mod(t)\right)e^{-C(\|x-\aa_1\|+\|x-\aa_2\|)}, 
 \eee
 where we defined the linear operator close to $V^{(N)}_j$ by: 
 \be
 \label{deflinear}
 L_j\Theta_j=-\Delta\Theta_j+\Theta_j+2\phi_{Re(\Theta_j\overline{V_j^{(N)}})}V_j^{(N)}+\left[\phi_{|V_j^{(N)}|^2}-\frac{\lambda_j^2\|Q\|_{L^2}^2}{4\pi\lambda_{j+1}\|\aa\|}\right]\Theta_j
 \ee
 Note that we have approximated the field created by the $V_{j+1}^{(N)}$ by the first order of its dipole expansion thanks to the space localization of $\Theta_j$  up to an $O(\frac{1}{\|\aa\|^2})$ term. We thus obtain the preliminary system of modulation equations:
 \bea
\label{systemmodu}
& & \frac{d}{dt} \left\{Im\int\e\overline{\Theta}\right\} =  Re\int\e\left[\overline{\frac{1}{\lambda_j^4}\left[i\lambda_j^2\partial_t\Theta_j-L_j\Theta_j\right](y_j)e^{-i\gamma_j(t)+i\beta_j(t)\cdot x}}\right]\\
\nonumber & - & 2Re\int\e\overline{\phi_{Re(\Theta \overline{R_j})}R_{j+1}}+\frac{1}{\lambda^3_j}Re(\int S_j^{(N)}\overline{\Theta_j}) +O\left(\frac{C_N}{\|\aa\|^{N+1}}+\frac{\|\e\|_{H^1}}{\|\aa\|^2}+C_N\|\e\|_{H^1}^2\right)
\eea

 {\bf Step 2.} Approximate null space. 
 
 If we replace $L_j$ by the exact linear operator close to the ground state, then the choice of orthogonality conditions (\ref{orthconditions}) corresponds to the null space of $L$ and generates modulation equations that are quadratic in $\e$. In our setting, we claim that the choice (\ref{orthconditions}) reproduces the null space structure for $L_j$ up to an $O(\frac{1}{\|\aa\|^2})$ which is enough for our analysis.
 
 Indeed, let us rewrite the equation of $V_j^{(N)}$ (\ref{eq:22}) as follows:
 \be
 \label{eqvjn}
  \Delta V^{(N)}_j-V^{(N)}_j-\left[\phi_{|V_j^{(N)}|^2}-\frac{\lambda_j^2\|Q\|_{L^2}^2}{4\pi\lambda_{j+1}\|\aa\|}\right]V^{(N)}_j=O\left(\frac{1}{\|\aa\|^2}\right)e^{-C\|x-\aa_j\|}.
  \ee
  We then let the group of symmetry of (\ref{hartree}) act on this equation and differentiate with respect to the symmetry parameters to get:
  \be
  \label{phase}
  L_j(iV_j^{(N)})= O\left(\frac{1}{\|\aa\|^2}\right)e^{-C\|x-\aa_j\|}, \  \ (phase),
  \ee
 \be
 \label{scaling}
 L_j(\Lambda V_j^{(N)})=2\left[-1+\frac{\lambda_j^2}{4\pi\lambda_{j+1}\|\aa\|}\right]V_j^{(N)}+O\left(\frac{1}{\|\aa\|^2}\right)e^{-C\|x-\aa_j\|}, \  \ (scaling),
 \ee
 \be
  \label{translation}
  L_j(\nabla V_j^{(N)})= O\left(\frac{1}{\|\aa\|^2}\right)e^{-C\|x-\aa_j\|}, \  \ (translation),
  \ee
 \be
  \label{Galilei}
  L_j(iy_jV_j^{(N)})= -2\nabla V_j^{(N)}+O\left(\frac{1}{\|\aa\|^2}\right)e^{-C\|x-\aa_j\|}, \  \ (Galilei).
  \ee
 We now compute the modulation equations from (\ref{orthconditions}) by using (\ref{systemmodu}) with successively $\Theta_j=iV_j^{(N)}, \Lambda V_j^{(N)}, \nabla V_j^{(N)}, iy_j V_j^{(N)}$. 
 Observe that for these four functions, 
 $$i\partial_t \Theta_j=O\left(\frac{Mod(t)}{\|\aa\|}+\frac{1}{\|\aa\|^2}\right)e^{-C\|x-\aa_j\|},$$ and thus (\ref{phase}), (\ref{scaling}), (\ref{translation}), (\ref{Galilei}) yield: 
\be
\label{ciwo}
Mod(t)\leq \sum_{j=1}^2\left|Re\int\e\overline{\phi_{Re(\Theta \overline{R_j})}R_{j+1}}\right|+\frac{CMod(t)}{\|\aa\|}+\frac{C_N}{\|\aa\|^{N+1}}+\frac{\|\e\|_{H^1}}{\|\aa\|^2}+C_N\|\e\|_{H^1}^2.
\ee
We now use the dipole expansion on the term $\phi_{Re(\Theta \overline{R_j})}$, as in section \ref{sec:1.2},  the orthogonality condition $Re(\int\e\overline{R_{j+1}})=0$ from (\ref{orthconditions}) 
and the nondegeneracy relations $\int \Lambda Q Q \neq 0$, $\int y_i Q \partial_i Q\neq 0$ to conclude:
$$Re\int\e\overline{\phi_{Re(\Theta \overline{R_j})}R_{j+1}}=\frac{C_j}{\|\aa\|}Re(\int\e\overline{R_{j+1}})+O\left(\frac{\|\e\|_{H^1}}{\|\aa\|^2}\right)=O\left(\frac{\|\e\|_{H^1}}{\|\aa\|^2}\right).$$ Injecting this into (\ref{ciwo}) yields:
\be
\label{estmod}
Mod(t)\leq \frac{C_N}{\|\aa\|^{N+1}}+\frac{\|\e\|_{H^1}}{\|\aa\|^2}+C_N\|\e\|_{H^1}^2.
\ee
  
{\bf Step 3.} Integration of the modulation equations.

We now prove the estimates on the parameters of (\ref{hypflowboost}) as a simple consequence of (\ref{estmod}).
We inject (\ref{hypflowboot}), (\ref{estintermhyp}), (\ref{estintermpar}) into (\ref{estmod}) to get for $T_*\geq T_0(N)$ large enough, 
\be
\label{controlofmodulation}
\forall t\in [T_*,T_n], \ \  Mod(t)\leq \frac{1}{t^{\frac{N}{4}}}.
\ee
From (\ref{computmb}), the definition of $\lambda_j^{(N)}$, $\beta^{(N)}$  (see Lemma \ref{defdynamic}) and  the explicit structure of $M_j^{(N)}, B_j^{(N)}$ of degree $\geq 2$, there holds:
\bee &&
|\ld_j-\ld_j^{(N)}|+\|\bd_j-\bd^{(N)}_j\| \\
&& \leq   \frac{1}{t^{\frac{N}{4}}}+|M_j^{(N)}(P)-M_j^{(N)}(P^{(N)})|+\|B_j^{(N)}(P)-B_j^{(N)}(P^{(N)})\|\\
&  & \leq \frac{1}{t^{\frac{N}{4}}}+ \frac{C}{\|\aa^{(N)}\|^3}\|\aa-\aa^{(N)}\|+\frac{C}{\|\aa^{(N)}\|^2}\left[\|\beta-\beta^{(N)}\|+\sum_{k=1}^2|\lambda_k-\lambda_k^{(N)}|\right]
  \leq   \frac{C}{t^{2+\frac{N}{8}}}.
\eee
Hence from (\ref{initialization}): $\forall t\in [T_*,T_n]$,
\be
\label{estsaclingbeta}
|\lambda_j(t)-\lambda_j^{(N)}(t)|+\|\beta_j(t)-\beta_j^{(N)}(t)\|\leq\frac{C}{Nt^{1+\frac{N}{8}}}.
\ee
We now integrate the translation parameter and find from (\ref{controlofmodulation}) and (\ref{estsaclingbeta}):
$$\|\aad_j-\aad_j^{(N)}\|\leq C \|\beta_j(t)-\beta_j^{(N)}(t)\|+\frac{1}{t^{\frac{N}{4}}}\leq\frac{C}{Nt^{1+\frac{N}{8}}}$$ from which:  
$\forall t\in [T_*,T_n]$,
\be
\label{contrtranslation}
\|\aa_j(t)-\aa_j^{(N)}(t)\|\leq \frac{C}{N^2t^{\frac{N}{8}}}.
\ee
Next, we estimate the phase parameters from (\ref{controlofmodulation}), (\ref{estsaclingbeta}) and (\ref{contrtranslation}):
\bee
|\gd_j-\gd^{(N)}_j| & \leq  & C\left(|\lambda_j-\lambda_j^{(N)}|+\|\beta_j-\beta_j^{(N)}\|+t\|\bd_j-\bd_j^{(N)}\|+\frac{\|\aa_j-\aa_j^{(N)}\|}{\|\aa^{(N)}\|^2}\right)+\frac 1 {t^{\frac N4}}\\
& \leq & \frac{C}{Nt^{1+\frac{N}{8}}}
\eee
and thus: $\forall t\in [T_*,T_n]$,
\be
\label{estpahseeiuhuiop}
|\gamma_j(t)-\gamma_j^{(N)}(t)|\leq \frac{C}{N^2t^{\frac{N}{8}}}.
\ee
Finally, (\ref{estsaclingbeta}), (\ref{contrtranslation}) and (\ref{estpahseeiuhuiop}) now yield the estimate for the parameters in (\ref{hypflowboost}) for $N$ large enough.


\subsection{Control of $\e$ in $H^1$}


We now focus onto the proof of the bootstrap estimate (\ref{hypflowboost}) on $\e$.\\
Let us consider the energy conservation: 
\bee
2{\cal{H}}_0 & = & 2{\cal{H}}(R+\e)=2{\cal{H}}(R)-2Re(\e,\overline{\Delta R-\phi_{|R|^2}R})\\
& + &\int|\nabla \e|^2+\int\phi_{|R|^2}|\e|^2-2\int|\nabla \phi_{Re(\e\Rb)}|^2+2\int\phi_{Re(\e\Rb)}|\e|^2-\frac{1}{2}\int|\nabla\phi_{|\e|^2}|^2.
\eee
We add to the quadratic and higher order terms of this conservation law a localization  of the two other conservation laws, namely $L^2$ norm and momentum conservation,  on each solitary wave. We obtain a functional which is almost conserved by the linear flow governing the equation for $\e$ (\ref{eqe}) (see \cite{MMnls} and \cite{MMTnls} for similar technique). 

Let $\Psi(x)=\Psi(x_{(1)})$ (recall the notation $x=(x_{(1)},x_{(2)},x_{(3)})$) denote a smooth nonnegative cutoff function such that $\Psi(x)=1$ if $x_{(1)}\leq -1$ and $\Psi(x)=0$ if $x_{(1)}\geq 1$.  We distinguish   the two cases of the asympotic dynamics $P^{\infty}$ being a hyperbola or a parabola:
\begin{enumerate}
\item Hyperbolic regime:  According to (\ref{estintermhyp}), we may find a constant $C=C_{P^{\infty}}$ such that 
$$\psi_1(t,x)=\Psi(\frac{x}{Ct}),\ \ \psi_2(t,x)=1-\psi_1(t,x)$$ localizes around each center of mass $\aa_j(t)$, that is $\psi_j(x)=1$ around $\aa_j(t)$ (i.e. in some region of the form $|x_{(1)}-(\aa_j)_{(1)}(t)|\leq C t$). We then let $$\zeta_j(t,x)=\psi_j(t,x).$$
\item Parabolic regime:  According to (\ref{estintermpar}) , we may find a constant $C=C_{P^{\infty}}$ such that 
$$\psi_1(t,x)=\Psi(\frac{x}{Ct^{\frac{2}{3}}}),\ \ \psi_2(t,x)=1-\psi_1(t,x)$$ localizes around $\aa_j(t)$ ($|x_{(1)}-(\aa_j)_{(1)}(t)|\leq C t^{\frac 23}$). In this case, we let $$\zeta_j(t,x)=\frac{1}{2}.$$
\end{enumerate}
The point is that the common notation $\zeta_j$ will allow us to perform part of the proof for the two regimes at once.

\medskip

In both cases, we define a functional in $\e(t)$ which aims at localizing in space the global conservation laws:
\bea
\label{defG}
\nonumber {\cal{G}}(\e(t)) & = & \int|\nabla \e|^2+\int\phi_{|R|^2}|\e|^2-2\int|\nabla \phi_{Re(\e\Rb)}|^2+2\int\phi_{Re(\e \Rb)}|\e|^2-\frac{1}{2}\int|\nabla\phi_{|\e|^2}|^2\\
& + & \sum_{j=1}^2\left[\left(\frac{1}{\lambda_j^2}+\|\beta_j\|^2\right)\int\zeta_j|\e|^2-2\beta_j\int \psi_j\cdot Im(\nabla\e\overline{\e})\right].
\eea
Note that for the parabolic case, the $L^2$ term does not contain a cut-off (the proof would not work, this is where we will use $\lambda^{\infty}_1=\lambda^{\infty}_2$), while the momemtum term does contain a cut-off. For the hyperbolic case, both term contain a cut-off so that any $\lambda_j^{\infty}$ is possible.

\medskip

We first claim the following coercivity property of the linearized energy which is a standard consequence of the variational characterization of the ground state and the choice of the orthogonality conditions (\ref{orthconditions}), together with the nondegeneracy of the kernel of the linearized operator $L$ close to $Q$ as proved by Lenzmann \cite{Lenzman}: 

\begin{Lemma}[Coercivity of the linearized energy]
\label{coercivityquadra}
There exists a constant $c_0>0$ such that under the orthogonality conditions (\ref{orthconditions}),  for any $t\geq T_0$ large enough,
\be
\label{estH}
{\cal{G}}(\e(t))\geq c_0\|\e(t)\|^2_{H^1}.
\ee
\end{Lemma}

{\bf Proof of Lemma \ref{coercivityquadra}.} See Appendix B.\medskip

We now claim the following energy estimates on $\e$ which is the core of the proof:

\begin{Lemma}[Global energy estimate on $\e$]
\label{energyestimate}
There holds for $N$, $T_*\geq T_0$ large enough:  $\forall t\in [T_*,T_n]$, 
\be
\label{estenergye}
\left|\frac{d}{dt}{\cal{G}}(\e(t))\right|\leq C\left(\frac{C_{_N}}{t^{N+1}}+\frac{\|\e(t)\|^2_{H^1}}{t}\right)
\ee
with constant $C>0$ independent of $N$. 
\end{Lemma}  

\begin{remark}
\label{rkkey}
Observe that (\ref{estH}) and (\ref{estenergye}) formally imply a structure: 
\be
\label{esgfiu}
\left|\frac{d}{dt}\|\e\|_{H^1}^2\right|\leq C\left(\frac{C_{_N}}{t^{N+1}}+\frac{\|\e(t)\|^2_{H^1}}{t}\right).
\ee A fundamental difficulty here is the gain $\frac{1}{t}$ only in the control of the error terms $O(\frac{\|\e(t)\|^2_{H^1}}{t})$ mostly induced by the localization of the conservation laws. This size comes in the hyperbolic case from the distance $\|\alpha(t)\|\sim t$. In the parabolic case, the situation is much worse $\|\alpha(t)\|\sim t^{\frac{2}{3}}$ and hence we need to restrict to the symmetric case, see estimates (\ref{estzetakey}), (\ref{estkeypsitwo}), (\ref{estkeypsitwod}) below. Reintegrating (\ref{esgfiu}) backwards from infinity, we see that we may bootstrap an information $\|\e\|^2_{H^1}\lesssim \frac{1}{t^{\frac{N}{2}}}$ only if $N$ is large enough, hence the requirement of an approximate solution to a large enough order.
\end{remark}

Let us postpone the proof of Lemma \ref{energyestimate} and complete the proof of Lemma \ref{boot}.
\medskip

{\bf Proof of Lemma \ref{boot}.}

We may now conclude the proof of (\ref{hypflowboost}). Here the key is the fact that the constant $N$ in (\ref{estenergye}) may be fixed as large as we want. We integrate (\ref{estenergye}) in time so that from (\ref{initialization}) and (\ref{hypflowboot}): $\forall t\in [T_*,T_n]$, \ \ 
$$\left|{\cal{G}}(\e(t))\right|\leq \int_t^{T_n}\left(\frac{C_{N}}{\tau^{N+1}}+\frac{C}{\tau^{1+\frac{N}{2}}}\right)d\tau\leq \frac{C_{N}}{t^{N}}+\frac{C}{Nt^{\frac{N}{2}}}\leq \frac{C}{Nt^{\frac{N}{2}}}$$ for $t\geq T_*\geq T_0(N)$ large enough, where $C$ is independent of $N$. From (\ref{estH}) and (\ref{hypflowboost}), this implies: $${c_0}\|\e(t)\|_{H^1}^2\leq \frac{C}{Nt^{\frac{N}{2}}}$$ and (\ref{hypflowboost}) follows by fixing $N$ such that $$Nc_0>8C.$$ This concludes the proof of the bootstrap Lemma \ref{boot} assuming Lemma \ref{energyestimate}.

\medskip

We now turn to the proof of Lemma \ref{energyestimate}. We point out the key point of gaining the $\frac{\|\e\|_{H^1}^2}{t}$ degeneracy in the RHS of (\ref{estenergye}), the power $t$ being sharply enough to close the bootstrap estimate as in the proof of Lemma \ref{boot}.\medskip

{\bf Proof of Lemma \ref{energyestimate}.} 
Let us split: $${\cal{G}}={\cal{G}}_1+{\cal{G}}_2+{\cal{G}}_3$$ with $${\cal{G}}_1= \int|\nabla \e|^2+\int\phi_{|R|^2}|\e|^2-2\int|\nabla \phi_{Re(\e \Rb)}|^2+2\int\phi_{Re(\e \Rb)}|\e|^2-\frac{1}{2}\int|\nabla\phi_{|\e|^2}|^2,$$$${\cal{G}}_2=\sum_{j=1}^2\left(\frac{1}{\lambda_j^2}+\|\beta_j\|^2\right)\int \zeta_j|\e|^2,\qquad {\cal{G}}_3=-2\int \beta_j\psi_j\cdot Im(\nabla\e\overline{\e}).$$

{\bf Step 1.} Estimate on  ${\cal{G}}_1$. 
We compute: 
\bee
\frac{d}{dt}{\cal{G}}_1(t) & = & -2Im\int i\partial_t\e\left[\overline{\Delta\e}-\phi_{|R|^2}\overline{\e}-2\phi_{Re(\e\overline{R})}\Rb-2\phi_{Re(\e \overline{R})}\overline{\e}-\phi_{|\e|^2}\Rb-\phi_{|\e|^2}\overline{\e}\right]\\
& + & 2\int|\e|^2\phi_{Re(\partial_tR\overline{R})}+4Re\int \e\partial_t\Rb\phi_{Re(\e\Rb)}+2Re\int\e\partial_t \Rb\phi_{|\e|^2}.
\eee
We now observe from (\ref{defR}) and from computations similar to (\ref{changevariables}) that: $\forall j=1,2$,
\bee
\partial_t R_j(t,x) & = & \frac{1}{\lambda_j^4}\left[-\lambda_j^2\partial_tV^{(N)}_j-\lambda_j\ld_j\Lambda V_j^{(N)} +i\lambda_j^3(\bd_j\cdot y_j)V^{(N)}_j\right .\\
& - & \left .\lambda_j\aad_j\cdot\nabla V^{(N)}_j -i \lambda_j^2(\gd_j-\bd_j\cdot \aa_j)V_j^{(N)}\right](\frac{x-\aa_j(t)}{\lambda_j(t)})e^{i(t-\gamma_j(t))}e^{i\beta_j(t)\cdot x}
\eee
 and thus from (\ref{hypflowboot}), (\ref{computV}) and (\ref{estxhyp}), (\ref{estxparab}): $\forall j=1,2$,
 \bea
 \label{estpartialtrj}
 \nonumber \partial_t R_j & = &  -2\beta_j\cdot\nabla R_j+i\left(\frac{1}{\lambda_j^2}+\|\beta_j\|^2\right)R_j+O\left(\frac{1}{\|\aa\|^2}+\frac{1}{t^{\frac{N}{8}}}\right)e^{-C_N\|x-\aa_j\|}\\
 & = & -2\beta_j\cdot\nabla R_j+i\left(\frac{1}{\lambda_j^2}+\|\beta_j\|^2\right)R_j+O\left(\frac{1}{t}\right)e^{-C_N\|x-\aa_j\|}.
\eea
We now use the localization properties of $R_j$ around $\aa_j$, inject the equation for $\e$ (\ref{eqe}) and get using (\ref{estpartialtrj}), the estimate of $\Psi^{(N)}$ (\ref{sizeerror}), the estimate on the geometrical parameters (\ref{controlofmodulation}) and the bootstrap estimate (\ref{hypflowboot}) for cubic and higher terms in $\e$ to derive:
$$
\frac{d}{dt}{\cal{G}}_1= \sum_{j=1}^2{\cal K}_{1,j}+
O\left(\frac{C_N\|\e(t)\|_{H^1}}{t^{N+1}}+\frac{\|\e(t)\|^2_{H^1}}{t}\right).
$$ with 
\bea
\label{estgunj}
 {\cal{K}}_{1,j}& = &   -4\int|\e|^2\phi_{Re(\beta_j \cdot \nabla R_j\overline{R}_j)}\\
\nonumber & - &  4\left(\|\beta_j\|^2+\frac{1}{\lambda_j^2}\right)\int Im(\e\overline{R}_j)\phi_{Re(\e\Rb)}-8\int Re(\e\beta_j\cdot\nabla \overline{R}_j)\phi_{Re(\e\overline{R})}. \eea

{\bf Step 2.} Estimate on ${\cal{G}}_2$. 
We claim that:
\bea
\label{estgtwo}
\frac{d}{dt}{\cal{G}}_2(t) & = & 4\sum_{j=1}^2 \left(\frac{1}{\lambda_j^2}+\|\beta_j\|^2\right)\int Im(\e\overline{R_j})\phi_{Re(\e\overline{R})}\\
\nonumber & + & O\left(\frac{C_N\|\e(t)\|_{H^1}}{t^{N+1}}+\frac{\|\e(t)\|^2_{H^1}}{t}\right).
\eea

We argue differently in the hyperbolic and parabolic regimes:
\smallskip

\noindent 1. Hyperbolic case: From (\ref{estintermhyp}),  the distance beween the two solitary waves in the hyperbolic regime is lower bounded: $\|\aa(t)\|\geq Ct$, 
$C>0$, which allows the estimates:
\be
\label{estzetakey}
\|\nabla \zeta_j(t)\|+|\partial_t \zeta_j|\leq \frac{C}{t}.
\ee
Moreover,  $\|\dot\beta_j(t)\|+|\dot\lambda_j(t)|\leq \frac {C}{t^2}$ from the equations of motion. Thus,
from (\ref{eqe}) using also (\ref{hypflowboot}), (\ref{controlofmodulation}): 
\bea
\label{choss}
\nonumber& & \frac{d}{dt}\left[\left(\frac{1}{\lambda_j^2}+\|\beta_j\|^2\right)\int\zeta_j|\e|^2\right]= 2\left(\frac{1}{\lambda_j^2}+\|\beta_j\|^2\right)\int Im(i\partial_t\e\zeta_j\overline{\e})+O\left(\frac{\|\e(t)\|_{H^1}^2}{t}\right)\\
\nonumber & = & \left(\frac{1}{\lambda_j^2}+\|\beta_j\|^2\right)\left[2\int  Im(\nabla\zeta_j\cdot\nabla \e\overline{\e})+4\int Im(\e\zeta_j\overline{R})\phi_{Re(\e\overline{R})}\right]\\
& + & O\left(\frac{C_N\|\e(t)\|_{H^1}}{t^{N+1}}+\frac{\|\e(t)\|^2_{H^1}}{t}+\|\e(t)\|_{H^1}Mod(t)\right),
\eea
and (\ref{estgtwo}) then follows from (\ref{estmod}).

\smallskip

\noindent 2. Parabolic case: From (\ref{estintermpar}), $\|\aa(t)\|\sim t^{\frac{2}{3}}$ in this case, and thus $\|\dot\beta_j(t)\|+|\dot\lambda_j(t)|\leq \frac {C}{t^{\frac 43}}$. Since $\zeta\equiv \frac 12$, using (\ref{estmod}) and arguing as above, we get  
\bee
 \frac{d}{dt}{\cal{G}}_2(t) &=&  2\sum_{j=1}^2 \left(\frac{1}{\lambda_j^2}+\|\beta_j\|^2\right)\int Im(\e\overline{R})\phi_{Re(\e\overline{R})}+ O\left(\frac{C_N\|\e(t)\|_{H^1}}{t^{N+1}}+\frac{\|\e(t)\|^2_{H^1}}{t}\right)\\
 & =&  4\sum_{j=1}^2 \left(\frac{1}{\lambda_j^2}+\|\beta_j\|^2\right)\int Im(\e\overline{R_j})\phi_{Re(\e\overline{R})}\\
\nonumber & + & O\left(\left[|\lambda_1-\lambda_2|+|\|\beta_1\|-\|\beta_2\||\right]\|\e(t)\|_{H^1}^2\right)+O\left(\frac{C_N\|\e(t)\|_{H^1}}{t^{N+1}}+\frac{\|\e(t)\|^2_{H^1}}{t}\right)
\eee
We now claim that:
\be
\label{hdwioqhoq}
|\lambda_1-\lambda_2|+|\|\beta_1\|-\|\beta_2\||\leq \frac{C}{t}
\ee
and (\ref{estgtwo}) follows.

\medskip

Proof of (\ref{hdwioqhoq}): We estimate from (\ref{mod}), (\ref{controlofmodulation}) and $\lambda_j(T_n)=\lambda^{(N)}_j(T_n)$ from (\ref{hartreeTn}): $\forall t\in [T_0,T_n]$,
\bea
\label{vobheoheoh}
\nonumber  |\lambda_1(t)-\lambda_2(t)| & \leq &   |\lambda_1(t)-\lambda_1^{(N)}(t)|+ |\lambda^{(N)}_1(t)-\lambda^{(N)}_2(t)|+|\lambda_2(t)-\lambda^{(N)}_2(t)|\\
\nonumber   & \lesssim  & \int_t^{T_n}Mod(\tau)d\tau+ |\lambda^{(N)}_1(t)-\lambda^{(N)}_2(t)|\\
 & \lesssim & \frac{1}{t^{\frac{N}{8}}}+|\lambda^{(N)}_1(t)-\lambda^{(N)}_2(t)|.
 \eea
We now estimate from the definition of $P_N$ and (\ref{estmjdeux}): 
\be
\label{fnohoif}
 \nonumber |\dot{\lambda}^{(N)}_1-\dot{\lambda}^{(N)}_2|  \lesssim  \left[\frac{ \|\beta^{(N)}\| }{\|\aa^{(N)}\|^2} \, |\lambda^{(N)}_1-\lambda^{(N)}_2|+\frac{1}{\|\aa^{(N)}\|^3}\right].
 \ee 
Now observe from Lemma \ref{defdynamic} and the symmetry assumption on the asymptotic trajectory in the parabolic case that
$$\lim_{t\to +\infty}\lambda_1^{(N)}(t)=\lambda_1^{\infty}=\lambda_2^{\infty}=\lim_{t\to +\infty}\lambda_2^{(N)}(t)$$ and thus the time integration of (\ref{fnohoif}) with
$\|\aa^{(N)}(t)\|\sim t^{\frac{2}{3}}$, $\|\beta^{(N)}\|\sim t^{-\frac 13}$   yield a first estimate $|\lambda_1^{(N)}(t)-\lambda_2^{(N)}(t)|\lesssim \frac{1}{t^{\frac{2}{3}}}$. Reinjecting this bound into (\ref{fnohoif}) we eventually derive: 
\be
\label{esyhiofhs}
|\lambda_1^{(N)}(t)-\lambda_2^{(N)}(t)|\lesssim \frac{1}{t}
\ee 
which together with (\ref{vobheoheoh}) yields (\ref{hdwioqhoq}) for $|\lambda_1-\lambda_2|$.

We argue similarily for $\beta$. From (\ref{mod}), (\ref{controlofmodulation}) and $\beta_j(T_n)=\beta^{(N)}_j(T_n)$, there holds: $\forall t\in [T_0,T_n]$,
\bea
\label{defdynamicpouetpouet}
\nonumber  \|\beta_1(t)+\beta_2(t)\| & \leq &   \|\beta_1(t)-\beta_1^{(N)}(t)\|+ \|\beta^{(N)}_1(t)+\beta^{(N)}_2(t)\|+\|\beta_2(t)-\beta^{(N)}_2(t)\|\\
& \lesssim &  \frac{1}{t^{\frac{N}{8}}}+\|\beta^{(N)}_1(t)+\beta^{(N)}_2(t)\|.
 \eea
We now estimate from the definition of $P_N$ and (\ref{modulationorder2}): 
$$
 \|\dot{\beta}^{(N)}_1+\dot{\beta}^{(N)}_2\|  \lesssim  \left[\frac{ |\lambda^{(N)}_1-\lambda^{(N)}_2|}{\|\aa^{(N)}\|^2}+\frac{1}{\|\aa^{(N)}\|^3}\right]\lesssim \frac{1}{t^{2}}
 $$
 where we used (\ref{esyhiofhs}) in the last step. Integrating this with $\lim_{t\to+\infty}\beta_1^{(N)}=\lim_{t\to+\infty}\beta_1^{(N)}=0$ from Lemma \ref{defdynamic}  yields: $$\|\beta_1^{(N)}+\beta_2^{(N)}\|\lesssim \frac{1}{t}$$ which together with (\ref{defdynamicpouetpouet}) concludes the proof of (\ref{hdwioqhoq}).

\medskip

{\bf Step 3.} Estimate on ${\cal{G}}_3$. 
First we have, since $\|\dot \beta(t)\|\leq C t^{-\frac 43}$,
\bea
\label{wqpw}
 \frac{d}{dt}\int \left[-2Im(\beta_j\psi_j\cdot \nabla\e\overline{\e})\right]&=&-4\int Re\left(i\partial_t\e\left[\frac{1}{2}\beta_j\nabla\psi_j\overline{\e}+\beta_j\psi_j\cdot\nabla\overline{\e}\right]\right)
 \nonumber \\ &-&\int 2Im(\beta_j\partial_t \psi_j\cdot \nabla\e\overline{\e})
 + O\left(\frac{\|\e(t)\|^2_{H^1}}{t}\right).\nonumber 
 \eea
 We now observe in the  hyperbolic case that (\ref{estintermhyp}) implies:
\be
\label{estkeypsitwo}
\|\beta_j(t)\|\|\nabla \psi_j(t)\|_{L^{\infty}}\leq \frac{C}{t}.
\ee 
A fundamental observation is that the same estimate holds in the parabolic case. Indeed, from (\ref{estintermpar}): $$\|\beta_j(t)\|\leq \frac{C}{t^{\frac{N}{8}}}+\|\beta_j^{(N)}(t)\|\leq \frac{C}{t^{\frac{1}{3}}}$$  and hence 
\be
\label{estkeypsitwod}
\|\beta_j(t)\|\|\nabla \psi_j(t)\|_{L^{\infty}}\leq \frac{C}{t^{\frac{1}{3}+\frac{2}{3}}}\leq \frac{C}{t}.
\ee 
This means that the localization in space of the kinetic momentum is harmless in the parabolic case even though both solitons are not far apart thanks to the explicit decay in time of the Galilean parameter. Note that by the same argument
$$\|\beta_j(t)\|\|\partial_t \psi_j(t)\|_{L^{\infty}}\leq \frac{C}{t}.$$

From the previous observations and  (\ref{eqe}), (\ref{hypflowboot}), (\ref{controlofmodulation}), we  derive
\bee
& & \frac{d}{dt}\left[-2Im(\beta_j\psi_j\cdot \nabla\e\overline{\e})\right]\\
&=&-4\int Re\left(i\partial_t\e\left[\frac{1}{2}\nabla\cdot(\beta_j\psi_j)\overline{\e}+\beta_j\psi_j\cdot\nabla\overline{\e}\right]\right)+O\left(\frac{\|\e(t)\|^2_{H^1}}{t}\right)\\
& = & -4\int Re\left([2\phi_{Re(\e\overline{R})}R+\phi_{|R|^2}\e]\beta_j\psi_j\cdot\nabla\overline{\e}\right)+O\left(\frac{C_N\|\e(t)\|_{H^1}}{t^{N+1}}+\frac{\|\e(t)\|^2_{H^1}}{t}\right)\\
& = & 8\int Re(\e\beta_j\psi_j{\cdot}\nabla\overline{R})\phi_{Re(\e\overline{R})}+4\int|\e|^2\phi_{Re( \psi_j\beta_j \cdot \nabla R\overline{R})}+O\left(\frac{C_N\|\e(t)\|_{H^1}}{t^{N+1}}{+}\frac{\|\e(t)\|^2_{H^1}}{t}\right)
\eee
and thus from the localization properties of $\psi_j$ and $R_j$:
\bea
\label{estgthree}
\frac{d}{dt}{\cal{G}}_3(t) & = & \sum_{j=1}^2 8\int Re(\e\beta_j\cdot\nabla\overline{R_j})\phi_{Re(\e\overline{R})}+4\int|\e|^2\phi_{Re(\beta_j \cdot \nabla R_j\overline{R_j})}\\
\nonumber & + & O\left(\frac{C_N\|\e(t)\|_{H^1}}{t^{N+1}}+\frac{\|\e(t)\|^2_{H^1}}{t}\right).
\eea

{\bf Step 4.} Conclusion.
We add up (\ref{estgunj}), (\ref{estgtwo}) and (\ref{estgthree}) to get exactly:
$$\frac{d}{dt}{\cal{G}}(t)=O\left(\frac{C_N\|\e(t)\|_{H^1}}{t^{N+1}}+\frac{\|\e(t)\|^2_{H^1}}{t}\right)$$ and 
(\ref{estenergye}) follows.\\
This concludes the proof of Lemma \ref{energyestimate}.


\appendix



\section{Proof of Lemma \ref{defdynamic}}


\subsection{Hyperbolic Dynamics}

Let $P^\infty(t)$ be chosen as at the beginning of section 3 in the hyperbolic regime.
Let $0<\epsilon<\frac 18$. Let $T_0>1$ large enough to be defined later.
We define 
$\Omega_{\epsilon,T_0}=\Omega$ the set of functions
$P(t)=(\aa_1(t),\aa_2(t),\beta_1(t),\beta_2(t),\lambda_1(t),\lambda_2(t))$
defined on  $[T_0,+\infty)$ such that 
$|||P|||\leq 1$
where the norm $|||.|||$ is defined as follows
\[
||| P |||= \sum_{j=1}^2\sup_{t\in [T_{0}, \infty)}\left\{t^{1-\epsilon}\|\aa_{j}(t)-\aa_{j}^\infty(t)\|+t^{2-\frac \epsilon 2}\|\beta_{j}(t)-\beta_j^\infty(t)\|+t^{1-\frac \epsilon 4}|\lambda_{j}(t)-\lambda_j^\infty|\right\}.
\]

Now, we consider the map $\Gamma$ which corresponds to solve the dynamical system (\ref{bofbof})
where the functions $B_{j}^{(N)}$, $M_j^{(N)}$ are defined in Proposition \ref{propinduction}.
For $P(t)\in \Omega$, we set
\[
\Gamma P(t)=(\Gamma \aa_{1} (t), \Gamma \aa_{2} (t),\Gamma \beta_{1} (t), \Gamma \beta_{2} (t),\Gamma\lambda_{1} (t),\Gamma\lambda_{2} (t)),\quad t\in [T_{0}, \infty),
\]
where  (as usual $\aa(t)=\aa_2(t)-\aa_1(t)$, $\beta(t)=\beta_2(t)-\beta_1(t)$)
\[
\Gamma \aa_{j}(t)=\aa_{j}^{\infty}(t)+\int_{t}^{\infty}2(-\beta_{j}(\tau)+\beta_j^{\infty}(\tau))d\tau,
\]
\[
\Gamma \beta_{j}(t)=\lim_{t\rightarrow\infty}\beta^{\infty}_{j}(t)-\int_{t}^{\infty}B_{j}^{(N)}(\aa(\tau), \beta(\tau),\lambda_{1}(\tau),\lambda_2(\tau))d\tau,
\]
\[
\Gamma\lambda_{j}(t)=\lambda^{\infty}_{j}-\int_{t}^{\infty}M_{j}^{(N)}(\aa(\tau),\beta(\tau),\lambda_{1}(\tau),\lambda_{2}(\tau))d\tau.
\]
Now we claim 

\begin{Lemma} There exists $T_{0}$ large enough such that $\Gamma$ maps $\Omega$ into $\Omega$. Moreover, we have the {\it{contraction property}}: let $P_{a}$, $P_{b}\in \Omega$, then
\[
|||\Gamma P_{a}-\Gamma P_{b}|||\leq \frac{1}{2} |||P_{a}- P_{b}|||.
\]
\end{Lemma}

We note that this lemma immediately implies Lemma~\ref{defdynamic}  via a standard Banach fixed point argument. 

\medskip

\textbf{Proof.} We show that $\Gamma$ maps $\Omega$ into $\Omega$. The contraction property follows in identical manner by forming the difference equations.

First, we consider $\Gamma \aa_{j}$. Here we have, by $|||P|||\leq 1$, 
\[
\|\Gamma \aa_{j}(t)-\aa_{j}^{\infty}(t)\|=\|\int_{t}^{\infty}2(\beta_{j}(\tau)-\beta_{j}^{\infty}(\tau))d\tau\|\lesssim t^{-1+\frac \e 2}.
\]
 Thus, choosing $T_{0}$ large enough, we get 
\[
\|\Gamma \aa_j(t)-\aa_{j}^{\infty}(t)\|\leq t^{-1+\epsilon}.
\]

Next, we consider $\Gamma\beta_{j}(t)$. Note from (\ref{defpinfty}) and (\ref{modulationorder2})  that we have 
\[
\beta^{\infty}_{j}(t)=\lim_{t\rightarrow\infty}\beta^{\infty}_{j}(t)-\int_{t}^{\infty}b_{j}^{(2)}( \aa^{\infty}(\tau),\beta^{\infty}(\tau),\lambda_{1}^{\infty},\lambda_{2}^{\infty})d\tau.
\]
Thus,
\bea\nonumber
&\Gamma\beta_{j}(t)-\beta^{\infty}_{j}(t)=\int_{t}^{\infty}[b_{j}^{(2)}( \aa^{\infty}(\tau),\beta^{\infty}(\tau),\lambda_{1}^{\infty},\lambda_{2}^{\infty})-b_{j}^{(2)}(\aa(\tau),\beta(\tau),\lambda_{1}(\tau),\lambda_{2}(\tau))]d\tau\\\nonumber
&-\sum_{n=3}^{N}\int_{t}^{\infty}b_{j}^{(n)}(\aa(\tau),\beta(\tau),\lambda_{1}(\tau),\lambda_{2}(\tau))d\tau.
\eea
Using the expression of $b_j^{(2)}$ in (\ref{modulationorder2}), the estimates on $b_j^{(n)}$ for $n\geq 3$ and the assumption $|||P(t)|||\leq 1$, we get (for $\epsilon>0$ small enough)
\[
\left\|\int_{t}^{\infty}[b_{j}^{(2)}(\aa^{\infty}(\tau),\beta^{\infty}(\tau),\lambda_{1}^{\infty},\lambda_{2}^{\infty})-b_{j}^{(2)}(\aa(\tau),\beta(\tau),\lambda_{1}(\tau),\lambda_{2}(\tau))]d\tau\right\|\lesssim t^{-2+\frac \e 4},
\]
\[
\left\|\sum_{n=3}^{N}\int_{t}^{\infty}b_{j}^{(n)}(\aa(\tau),\beta(\tau),\lambda_{1}(\tau),\lambda_{2}(\tau), )d\tau\right\|\lesssim t^{-2},
\]
whence we have shown 
$
\|\Gamma\beta_{j}(t)-\beta^{\infty}_{j}(t)\|\lesssim t^{-2+\frac \e 4}.
$
Choosing $T_{0}$ large enough, for $t\geq T_{0}$, we infer
\[
\|\Gamma\beta_{j}(t)-\beta^{\infty}_{j}(t)\|\leq t^{-2+\frac \e 2}
\]

Finally, we consider $\Gamma\lambda_{j}(t)$. Since $m_j^{(2)}\in {\cal S}_2$, we 
have $|m_j^{(2)}|\lesssim t^{-2}$, thus we get
$
|\Gamma\lambda_{j}(t)-\lambda_{j}^{\infty}|\lesssim t^{-1+\frac \e 4}
$
whence 
\[
|\Gamma\lambda_{j}(t)-\lambda_{j}^{\infty}(t)|\leq t^{-1+\epsilon}
\]
if $T_{0}$ is large enough. This concludes the proof.


\subsection{Parabolic Dynamics}


Here we prove the following: 

\begin{Lemma} [Existence of almost parabolic trajectories]
\label{lemmaappendixa2}  Let $\lambda_{1}^{\infty}=\lambda_{2}^{\infty}=\lambda^{\infty}$ and  $P^{\infty}(t)(\aa_1^{\infty}(t),\aa_2^{\infty}(t),\beta_1^{\infty}(t),\beta_2^{\infty}(t))$ be a solution of the two-body problem (\ref{defpinfty}) of parabolic type satisfying (\ref{cnscnoson}).  There exists a solution $(\aa_1(t),\aa_2(t),\beta_1(t),\beta_2(t),\lambda_1(t), \lambda_2(t))$ of  (\ref{bofbof}) defined for $t\geq T_{0}$ with $T_{0}$ large enough satisfying:
$\forall j=1,2$,  
\be
\label{convergenceparappendix}
\sum_{k=1}^2\left|\frac{(\aa_j^{\infty})_{(k)}(t)}{(\aa_j)_{(k)}(t)}-1\right|+|(\aa_j)_{(3)}(t)|+\|\beta_j(t)-\beta_j^{\infty}(t)\|+|\lambda_j(t)-\lambda^{\infty}|\leq \frac{C_{P^{\infty}}}{t^{\frac{1}{3}}}
\ee
and 
\be
\label{estgaliela}
|\lambda_1(t)\aa_2(t)+\lambda_2(t)\aa_1(t)|\to 0 \ \ \mbox{as} \ \ t\to +\infty.
\ee
In particular, the asymptotic trajectory is given by the same parabola like for $P^{\infty}$.
\end{Lemma}

\textbf{Proof.}
{\bf Step 1.} Leading order dynamics.

First, we claim that the equation of motion of $P^{(N)}$ (i.e. (\ref{bofbof})) can be rewritten as follows:
\be\label{43bis}
\left   \{ \begin{array}{lllll}
	\ld_j=\frac{c_3 (\aa\cdot\beta) }{\|\aa\|^3}+\frac{g_j(\aa,\beta,\lambda_1,\lambda_2)}{\|\aa\|^2}, \ \ j=1,2\\
         \aad_j=2\beta_j, \ \ 
         \bd_j=\frac{(-1)^j}{4}\left[-\frac{c_{0}\aa}{||\aa||^{3}}+\frac{c_1 \aa}{\|\aa\|^{4}}+\Sigma_{j=1}^2\frac{c_{2}\aa(\lambda_j-\lambda^{\infty})}{\|\aa\|^{3}}\right]+\frac{f_j(\aa, \beta, \lambda_1,\lambda_2)}{\|\aa\|^{3}},
         \end{array}
\right .
\ee
where $c_0, c_1,c_2, c_3$ are explicit functions of $(\lambda^{\infty},\aa,\beta)$. Note that  the same constants $c_0,c_1,c_2,c_3$ appear for both $j=1,2$, because of the assumption $\lambda_1^{\infty}=\lambda_2^{\infty}=\lambda^{\infty}$. Moreover, the functions  $f_j, g_j$ satisfy the bound 
\bea
\label{estf}
&&|f_j| \lesssim  {\|\aa\|^{-1}}+|\beta|+\Sigma_{j=1}^2\left[|\lambda_j-\lambda^{\infty}|+\|\aa\|(\lambda_j-\lambda^{\infty})^{2}\right],
 \\ && |g_j|\lesssim {\|\aa\|^{-1}} +\Sigma_{j=1}^2\left[|\lambda_j-\lambda^{\infty}|+\|\aa\|(\lambda_j-\lambda^{\infty})^{2}\right]
\eea
as long as $\|\aa\|\geq 1$, $|\lambda_j-\lambda^{\infty}|+\|\beta\|\leq 1$.

Indeed, the structure of (\ref{43bis}) is deduced from a expansion of 
  the gravitational field (\ref{fieldnk}) to the order $ {\|\alpha\|^{-3}}$, and from structural symmetry properties:   both $T_j^{(1)}$ and $Re(T_j^{(2)})_{j=1,2}$ are chosen radially symmetric (see Remark \ref{rk1}).

\medskip

Let us  introduce an intermediate system $P^{(app)}$:

$$P^{(app)} \left   \{ \begin{array}{ll}
	\ld^{(app)}_1=\ld^{(app)}_2=\ld^{(app)}=\frac{c_3(\aa^{(app)}.\beta^{(app)})}{\|\aa^{(app)}\|^3},\\
         \aad^{(app)}_j=2\beta^{(app)}_j, \ \ \bd^{(app)}_j=\frac{(-1)^j}{4}\left[-\frac{c_{0}\aa^{(app)}}{\|\aa^{(app)}\|^{3}}+\frac{c_1 \aa^{(app)}}{\|\aa^{(app)}\|^{4}}+\frac{2 c_{2}\aa^{(app)}(\lambda^{(app)}-\lambda^{\infty})}{\|\aa^{(app)}\|^{3}}\right],
                 \end{array}
\right .
$$
and which implies the law for $\aa^{(app)}$:
\be
\label{lawforxapp}
\aad^{(app)}=2\beta^{(app)}, \ \ \ddot{\aa}^{(app)}=-\frac{c_{0}\aa^{(app)}}{||\aa^{(app)}||^{3}}+\frac{c_{1}\aa^{(app)}}{||\aa^{(app)}||^{4}}+\frac{2c_{2}\aa^{(app)}(\lambda^{(app)}-\lambda^{\infty})}{||\aa^{(app)}||^{3}}.
\ee
We claim that there exists a solution to $P^{(app)}$ which asympotically converges to $P^{\infty}$ in the sense of (\ref{convergenceparappendix}).

\medskip

Let us first observe from (\ref{lawforxapp}) that $\frac{d}{dt}(\aa^{(app)}\wedge \aad^{(app)})=0$ and hence the movement is planar. We assume that 
$\aa^{(app)}$ stays in the plane $\{x_{(3)}=0\}$.
 Choosing polar coordinates $(r,\theta)$ in this plane, $\aa^{(app)}(t)=(r(t)\cos \theta(t), r(t)\sin \theta(t),0)$, we rewrite the system as $$ \left  \{ \begin{array}{lll}
	r^2\dot{\theta}=a_0>0,\\
	\ld^{(app)}=\frac{c_3  \dot r}{r^2},\\
	\ddot{r}-r \, \dot{\theta}^2=-\frac{c_{0}}{r^2}+\frac{c_1}{r^3}+\frac{2c_{2}(\lambda^{(app)}-\lambda^{\infty})}{r^2}.
                 \end{array}
\right .
$$
We now use the classical change of variables $u=\frac{1}{r}$ and rewrite the system as an ODE for $u=u(\theta)$. From $$\frac{d}{dt}=\frac{d\theta}{dt}\frac{d}{d\theta}=a_0u^2\frac{d}{d\theta},$$ we compute: $$\dot{r}=-a_0\frac{du}{d\theta},\ \ \ddot{r}=-a_0^2u^2\frac{d^2u}{d\theta^2},\ \ \dot{\lambda}_j=a_0u^2\frac{d\lambda_j}{d\theta},$$ and hence the new system: 
$$ \left  \{ \begin{array}{ll}
	\frac{d\lambda^{(app)}}{d\theta}=- \frac 12 {c_3}\frac {du}{d\theta},\\
	a_0^2(\frac{d^2u}{d\theta^2}+u)=c_{0}-c_1 u-2c_{2}(\lambda^{(app)}-\lambda^{\infty}).
                 \end{array}
\right .
$$
We are led to compare the solution $\aa^\infty$ of  
$$ \left  \{ \begin{array}{lll}
	(r^{\infty})^2\dot{\theta}^{\infty}=a_0\\
	a_0^2(\frac{d^2u^{\infty}}{d(\theta^{\infty})^2}+u^{\infty})=c_{0}\\
	u^{\infty}(0)=\frac{du^{\infty}}{d\theta^{\infty}}(0)=0.
                 \end{array}
\right . \ \ \mbox{i.e.} \ \ r^{\infty}(\theta^{\infty})=\frac{a_0^2}{c_0(1-\cos(\theta^{\infty}))}.$$ 
with  $\aa^{(app)}$   solution on $(-\theta_0,0]$ of the system: 
$$ \left  \{ \begin{array}{ll}
	\frac{d\lambda^{(app)}}{d\theta}=-\frac 12 {c_3}\frac {du}{d\theta}, \quad \lambda^{(app)}(0)=\lambda^{\infty},
	\quad \lambda^{(app)}(\theta)-\lambda^\infty=-\frac 12 c_3 u(\theta), \\
	a_0^2(\frac{d^2u}{d\theta^2}+u)=c_0-(c_1-2c_2c_3) u=c_0-c_4 u,\quad  	u(0)=\frac{du}{d\theta}(0)=0.
                 \end{array}
\right .
$$
In fact, the system reduces to a linear one and the solution is explicit, but we do not need analytic formulas. The solution admits the following asymptotics near $\theta=0$: $$u(\theta)=\frac{c_0}{2a_0^2}\theta^2(1+O(\theta)), \ \  r(\theta)=\frac{2a_0^2}{c_0\theta^2}(1+O(\theta)),\ \  \lambda^{(app)}(\theta)=\lambda^{\infty}-\frac{c_0c_3}{2a_0^2}\theta^2(1+O(\theta)).$$
We compute the asymptotic time dependence using $$a_0=r^2\dot{\theta}=\frac{4a_0^4\dot{\theta}}{c^2_0\theta^4}(1+O(\theta))$$ which yields: 
\be
\label{estrtheta}
\theta(t)=-\left(\frac{4a_0^3}{3c_0^2t}\right)^{\frac{1}{3}}(1+O(\frac{1}{t^{\frac{1}{3}}})), \ \ r(t)=\left(\frac{9c_0}{2}\right)^{\frac{1}{3}}t^{\frac{2}{3}}(1+O(\frac{1}{t^{\frac{1}{3}}})),
\ee
\be
\label{estrthetal}
\lambda^{(app)}(t)-\lambda^{\infty}=-c_3\left(\frac{9c_0}{2}\right)^{-\frac{1}{3}}t^{-\frac{2}{3}}(1+O(\frac{1}{t^{\frac{1}{3}}})).
\ee
Similarily, we compute: 
\be
\label{estrthetainfty}
\theta^{\infty}(t)=-\left(\frac{4a_0^3}{3c_0^2t}\right)^{\frac{1}{3}}(1+O(\frac{1}{t^{\frac{1}{3}}})), \ \ r^{\infty}(t)=\left(\frac{9c_0}{2}\right)^{\frac{1}{3}}t^{\frac{2}{3}}(1+O(\frac{1}{t^{\frac{1}{3}}})),
\ee
and hence the polar curves $$(r(t)\cos(\theta(t)), r(t)\sin(\theta(t)))\quad  \mbox{and}\quad  (r^{\infty}(t)\cos(\theta^{\infty}(t)), r^{\infty}(t)\sin(\theta^{\infty}(t)))$$ are asympotic to each other as $t\to +\infty$. In particular, we have the asymptotic behavior of the cartesian coordinates of $\aa^{(app)}$: 
\be
\label{cartesainxapp}
\aa^{(app)}_{(1)}(t)=\left(\frac{9c_0t^2}{2}\right)^{\frac{1}{3}}(1+O(\frac{1}{t^{\frac{1}{3}}})), \ \
  \aa^{(app)}_{(2)}(t)=-\left(\frac{6a_0^3t}{c_0}\right)^{\frac{1}{3}}(1+O(\frac{1}{t^{\frac{1}{3}}})), \ \ \aa^{(app)}_{(3)}(t)=0.
\ee
We now solve for $(\aa_1(t), \aa_2(t))$ solution to $P^{(app)}$ by letting 
\be
\label{fhsfhowi}
\aa_2=\frac{\aa}{2}, \ \ \aa_1=-\frac{\aa}{2}
\ee 
thanks to the symmetry in $j=1,2$ of the system $P^{(app)}$, and hence the $\aa^{(app)}(t)$ curve also has a center of mass fixed at the origin. Equations (\ref{estrtheta}), (\ref{estrthetal}),  (\ref{estrthetainfty}) now imply (\ref{convergenceparappendix}) for $P^{(app)}$.

\medskip

{\bf Step 2.} Solving the full system.

We now claim that we can find a solution to the full system $P^{(N)}$ such that 
\be
\label{convegrencextxapp}
\Sigma_{j=1}^2\|\aa_j(t)-\aa_j^{(app)}(t)\|\to 0 \ \ \mbox{as} \ \ t\to +\infty.
\ee
 Recall that $P^{(N)}$ has been rewritten so that the following holds 
 $$\left   \{ \begin{array}{ll}
	\ld_j=\frac{c_3(\aa\cdot \beta)}{\|\aa\|^3}+\frac{g_j(\aa,\beta,\lambda_1,\lambda_2)}{\|\aa\|^2}, \ \ j=1,2\\
         \aad=2\beta, \ \  2\bd=-\frac{c_{0}\aa}{||\aa||^{3}}+\frac{c_1 \aa}{\|\aa\|^{4}}+\Sigma_{j=1}^2\frac{c_{2}\aa (\lambda_j-\lambda^{\infty})}{\|\aa\|^{3}}+\frac{f(\aa, \beta, \lambda_1,\lambda_2)}{\|\aa\|^{3}}, \ \ j=1,2.\\
         \end{array}
\right .
$$
where $f=f_2-f_1$ satisfies (\ref{estf}). We let $\aa(t)=\aa^{(app)}(t)+\dde(t)$, $\lambda_j(t)=\lambda^{(app)}(t)+\mu_j(t)$ and linearize the system around $\aa^{(app)}, \lambda^{(app)}$. Let us show how to bootstrap an estimate 
\be
\label{estboot}
\|\dde(t)\|\leq \frac{K}{t^{\frac{1}{4}}}, \ \ \|\dot{\dde}(t)\|\leq \frac{K}{t^{\frac{5}{4}}}.
\ee
\be
\label{estbootmu}
|\mu_1(t)|+|\mu_2(t)|\leq \frac{K}{t}.
\ee
for some large enough $K$; the claim then follows by a fixed point argument as in Appendix A.1, and thus will be omitted.

\medskip

The main linear term is expanded as follows
 $$\frac{\aa}{\|\aa\|^3}=\frac{\aa^{(app)}+\dde}{\|\aa^{(app)}+\dde\|^3}
 =\frac{\aa^{(app)}}{\|\aa^{(app)}\|^3}+\frac{1}{\|\aa^{(app)}\|^3}\left[\dde-3\frac{\aa^{(app)}\cdot \dde}{\|\aa^{(app)}\|^2}\aa^{(app)}\right]+ \frac{\|\dde\|^2O(1)}{\|\aa^{(app)}\|^4}.$$ We now observe from the explicit law (\ref{cartesainxapp}) that $$\frac{\aa^{(app)}\cdot \dde}{\|\aa^{(app)}\|^2}\aa^{(app)}=\dde_{(1)}+\frac{\|\dde\|O(1)}{t^{\frac{2}{3}}}$$ and so
 $$-\frac{c_0\aa}{\|\aa\|^3} -\frac{c_0\aa^{(app)}}{\|\aa^{(app)}\|^3} +\frac 1 {9t^2} \left(\begin{array}{lll}
	4\dde_{(1)}\\
	-2\dde_{(2)}\\
	-2\dde_{(3)}
	 \end{array}\right)
 +\frac{1}{t^{2+\frac{1}{3}}}O(\|\dde\|+\|\dde\|^2)
$$
Using the bootstrap (\ref{estboot}), (\ref{estf}) and (\ref{cartesainxapp}), the linearized system takes the form: 
 \be
 \label{estecayf}
 \left\{\begin{array}{ll}
 \dot{\mu}_j=\displaystyle\frac 1{t^{\frac 53}}O(1), \ \ j=1,2,\\
  \ddot{\dde}= \displaystyle\frac 1 {9t^2} \left(\begin{array}{lll}
 {4\dde_{(1)}}\\
	-{2\dde_{(2)}}\\
	-{2\dde_{(3)}}
	 \end{array}
\right) +F(t) 
\end{array}
\right .
\ \ \mbox{with} \ \ F(t)=\frac{1}{t^{2+\frac{1}{3}}}O(1).
\ee
An explicit solution is given by: $$ \dde(t)=\left(\begin{array}{lll}
 -\frac{3}{5}\left[t^{\frac{4}{3}}\int_{t}^{+\infty}\frac{F_{(1)}(\tau)}{\tau^{\frac{1}{3}}}d\tau+\frac{1}{t^{\frac{1}{3}}}\int_0^t\tau^{\frac{4}{3}}F_{(1)}(\tau)d\tau\right]\\
	-\frac{1}{3}\left[t^{\frac{2}{3}}\int_{\tau}^{+\infty}F_{(2)}(\tau)\tau^{\frac{1}{3}}d\tau-t^{\frac{1}{3}}\int_{\tau}^{+\infty}F_{(2)}(\tau)\tau^{\frac{2}{3}}d\tau\right]\\
	-\frac{1}{3}\left[t^{\frac{2}{3}}\int_{\tau}^{+\infty}F_{(3)}(\tau)\tau^{\frac{1}{3}}d\tau-t^{\frac{1}{3}}\int_{\tau}^{+\infty}F_{(3)}(\tau)\tau^{\frac{2}{3}}d\tau\right]
	 \end{array}
\right )
$$
from which (\ref{estboot}) is obtained with better constant.

To prove (\ref{estbootmu}), we have from $P^{(N)}$: 
\bee
\sum_{j=1}^2 |\dot{\lambda}_j-\dot{\lambda}^{(app)}| & \lesssim & \frac{1}{\|\aa^{(app)}\|^2}(\frac{\|\aa-\aa^{(app)}\|}{\|\aa^{app}\|}+\sum_{j=1}^2 |\lambda_j-\lambda^{(app)}|)+O(\frac{1}{\|\aa^{(app)}\|^3})\\
& \lesssim & \frac{1}{t^{\frac{4}{3}}}\Sigma_{j=1}^2 |\lambda_j-\lambda^{(app)}|)+\frac{1}{t^2}O(1)
\eee
which implies (\ref{estbootmu}) by integration in time.\\
We may now prove (\ref{estgaliela}). Indeed, 
\bee
\|\lambda_2(t)\aa_1(t)+\lambda_1(t)\aa_2(t)|\lesssim (\|\mu_2(t)\|+\mu_1(t)\|)t^{\frac{2}{3}}+\|\dde(t)\|\lesssim \frac{1}{t^{\frac{1}{4}}}
\eee
where we used $\aa_2^{(app)}+\aa_1^{(app)}=0$ from (\ref{fhsfhowi}).

This concludes the proof of Lemma \ref{lemmaappendixa2}.



\section{Proof of Lemma \ref{coercivityquadra}}

The proof follows exactly the same lines as the proof of Lemma 2.6 in \cite{MMTnls}.
We proceed into four steps. 

\medskip

\textbf{Step 1.} First, we claim that there exists $c>0$ such that for any real-valued $v\in H^1$
\be\label{surL+}
	(v,Q)=(v,yQ)=0 \quad \Rightarrow \quad (L_+ v,v)\geq c \|v\|_{H^1}^2,
\ee
\be\label{surL-}
	(v,\Lambda Q)=(v,\nabla Q)=0 \quad \Rightarrow \quad (L_- v,v)\geq c \|v\|_{H^1}^2.
\ee
Indeed, estimates (\ref{surL+})-(\ref{surL-}) follow from the arguments of \cite{W0},
proof of Propositions 2.9 and 2.10 and the fact that the kernels of $L_\pm$ are explicitely known from \cite{Lenzman} (the kernel of $L_-$ is spanned by $Q$ and the kernel of $L_+$ is spanned by the components of $\nabla Q$).

\medskip

\textbf{Step 2.} 
Second, we claim that for all  $\e\in H^1$ satisfying the orthogonality conditions (\ref{defej})-(\ref{orthconditions}), and $T_0$ large enough, we have
\be
\label{Gj}
{\cal{G}}_j(\e)\geq c_0\|\e\|^2_{H^1},
\ee
where
\bea
\label{defGjj}
\nonumber {\cal{G}}_j(\e) & = & \int|\nabla \e|^2+\int\phi_{|R_j|^2}|\e|^2-2\int|\nabla \phi_{Re(\e \bar R_j)}|^2\\
& + &   \left(\frac{1}{\lambda_j^2}+\|\beta_j\|^2\right)\int  |\e|^2-2\int  \beta_j\cdot Im(\nabla\e\overline{\e}) .
\eea
Let 
$Q_j(t,x)= \frac{1}{\lambda_j^2(t)}Q(\frac{x-\aa_j(t)}{\lambda_j(t)})e^{-i\gamma_j(t)+i\beta_j(t)\cdot x}$.
By the definition of $V_j^{(N)}$ in Proposition \ref{propinduction},  the definition of $R_j$
in Proposition \ref{constructionansatz}, and the properties of the functions $T_j^{(n)}$
we have 
$$
|(R_j- Q_j)(t,x)|\leq \frac C {\|\aa\|} e^{-\gamma \|x\|} \leq \frac C {t^{\frac 23}} e^{-\gamma \|x\|}. 
$$
Thus, if we define
\bea
\label{defGjjt}
\nonumber \tilde {\cal{G}}_j(\e) & = & \int|\nabla \e|^2+\int\phi_{|Q_j|^2}|\e|^2-2\int|\nabla \phi_{Re(\e \bar Q_j)}|^2\\
& + &   \left(\frac{1}{\lambda_j^2}+\|\beta_j\|^2\right)\int  |\e|^2-2\int  \beta_j\cdot Im(\nabla\e\overline{\e}) ,
\eea
we have
$$|\tilde {\cal{G}}_j(\e)-{  \cal{G}}_j(\e)|\leq \frac C {t^{\frac 23}} \|\e\|_{H^1}^2.$$

 Then, by standard computations (see e.g. proof of Claim 7 in \cite{MMTnls}), we have
$$
\tilde {\cal{G}}_j(\e) = (L_+ Re(\e_j), Re(\e_j))+ (L_- Im(\e_j),Im(\e_j)) + O(\frac 1{\|\aa(t)\|}) \|\e\|_{H^1}^2,
$$
where $\e_j$ are defined in (\ref{defej}).
Thus,  (\ref{Gj}) follows from (\ref{surL+})-(\ref{surL-}) and
$\|\e \|_{H^1} \leq C \|\e_j\|_{H^1}$.
 
\medskip

\textbf{Step 3.} The function $\psi_j$ being defined as in Section \ref{sec:3.2}, let
\bea
\label{defGj}
\nonumber {\cal{G}}_{j,\psi}(\e) & = & \int \psi_{j} |\nabla \e|^2+\int\phi_{|R_j|^2}|\e|^2-2\int|\nabla \phi_{Re(\e \bar R_j)}|^2\\
& + &  \left[\left(\frac{1}{\lambda_j^2}+\|\beta_j\|^2\right)\int \psi_{j} |\e|^2-2\int  \psi_{j}\beta_j\cdot Im(\nabla\e\overline{\e})\right].
\eea
Then, we claim that
\be\label{eeb}
\nonumber {\cal{G}}_{j,psi}(\e)\geq c \int \psi_{j} (|\nabla \e|^2+|\e|^2)-\frac C {t^{\frac 23}} \|\e\|_{H^1}^2.
\ee
Proof of (\ref{eeb}).
Let $\tilde \e_{j}=\e \sqrt{\psi_{j}}$.
Proceeding as in the proof of Claim 8 in \cite{MMTnls}, we have
$$\int |\nabla \tilde \e_{j}|^2 - \frac C {t^{\frac 23}} \|\e\|_{H^1}^2
\leq \int  \psi_{j}|\nabla \e|^2 
\leq \int |\nabla \tilde \e_{j}|^2 + \frac C {t^{\frac 23}} \|\e\|_{H^1}^2$$
and $\psi_{j}\beta_j\cdot Im(\nabla\e\overline{\e})
= \beta_j\cdot Im(\nabla{\tilde \e_{j}}\overline{\tilde \e_{j}})$.

Next, we have
$$\int  \phi_{|R_j|^2}|\e|^2 = \int  \phi_{|R_j|^2}|\tilde \e_{j}|^2
+ \int \phi_{|R_j|^2}|\e|^2 (1-\psi_{j})$$
and
\bea
- \int|\nabla \phi_{Re(\e \bar R_j)}|^2
&=&\int \phi_{Re(\e \bar R_j)} {Re(\e \bar R_j)}
= \int \phi_{Re(\tilde \e_{j} \bar R_j)} {Re(\e \bar R_j)}\\
&+& \int \phi_{Re(\e (1-\psi_{j})\bar R_j )} {Re(\e(1+\psi_{j})\bar R_j)}.
\eea
Thus, by the decay properties of $\psi$ and $R_j$, we obtain
\bee
\left|\int\phi_{|R_j|^2}|\e|^2-2\int|\nabla \phi_{Re(\e \bar R_j)}|^2
- \int\phi_{|R_j|^2}|\tilde \e_{j}|^2+2\int|\nabla \phi_{Re(\tilde \e_{j}\bar R_j)}|^2\right|
 \leq \frac C{t^{\frac 23}} \|\e\|_{H^1}^2.
\eee
From these estimates, we obtain 
\be\label{eeeeb}
\nonumber {\cal{G}}_{j,\psi}(\e)\geq {\cal{G}}_{j}(\tilde \e_{j})
- \frac C{t^{\frac 23}} \|\e\|_{H^1}^2.\ee
By standard arguments (see e.g. \cite{MMnls}), since for $t$ large $\tilde \e_{j}$ almost  satisfies   the same orthogonality conditions as $\e$, we obtain by step 2,
${\cal{G}}_{j}(\tilde \e_{j})\geq c \|\tilde \e_{j}\|_{H^1}^2 
- \frac C{t^{\frac 23}} \|\e\|_{H^1}^2$ and (\ref{eeb}) follows.
\medskip

\textbf{Step 4.} Conclusion. We decompose ${\cal{G}}(\e)$ as follows
\bea
\label{defGjjj}\nonumber
{\cal{G}}(\e) & = & \sum_{j=1}^2 
{\cal{G}}_{j,\psi}(\e)
\\  \label{g2}
&+& \int \phi_{|R|^2-(|R_1|^2|+|R_2|^2)} |\e|^2 \\ \label{g4} &-&
2 \int (|\nabla \phi_{Re(\e \bar  R)} |^2 -|\nabla \phi_{Re(\e \bar R_1)} |^2
-|\nabla \phi_{Re(\e\bar R_2)} |^2) \\
&+&\label{g3}
2 \int \phi_{Re(\e \bar R)} |\e|^2 - \frac 12 \int |\nabla \phi_{|\e|^2}|^2.
\eea
By the decay properties of $R_1$, $R_2$, the term in (\ref{g2}) is easily controlled by
$C/t$. The terms in (\ref{g3}) are nonlinear in $\e$, and we obtain by the Hardy-Littlewood-Sobolev inequality
$$
\left|  \int \phi_{Re(\e R)} |\e|^2 \right|+ \frac 12 \int |\nabla \phi_{|\e|^2}|^2
\leq C (\|\e\|_{H^1}^3+\|\e\|_{H^1}^4) 
\leq \frac C{t^{\frac N4}}\|\e\|_{H^1}^2,
$$
where we have used (\ref{hypflowboot}) for the last estimate.

Now, we estimate (\ref{g4}). 
First, by the decay properties of $R_j$, we have 
$$\left|\int \phi_{Re(\e(t) \bar R_1)} Re(\e(t)\bar R_2(t))\right|\leq
C \int e^{- \gamma |y-\aa_1(t)|} |\e(t,y)| e^{- \gamma |z-\aa_2(t)|} |\e(t,z)|\frac 1{|y-z|} dy dz.
$$
We cut this integral into two pieces depending on  $|y-z\|\geq \frac 12 \|x\|$ or $\|y-z\|\leq \frac 12 \|x\|$.
The first part, i.e. when $\|y-z\|\geq \frac 12 \|x\|$, is easily estimated by
$\frac C{\|\aa\|} \|\e\|_{L^2}^2$ using Cauchy-Schwarz inequality.
For the second part, we observe that if $\|y-z\|\leq \frac 12 \|\aa\|$
then $\|z-\aa_2\|\geq \|\aa_2-\aa_1\|-\|y-\aa_1\|-\|z-y\|\geq \frac 12 \|\aa\| - \|y-\aa_1\|$ and so
$$ e^{- \gamma \|y-\aa_1(t)\|}e^{- \gamma |z-\aa_2(t)|}\leq e^{- \gamma \|y-\aa_1(t)\|} e^{- \frac 12 \gamma \|z-\aa_2(t)\|}
\leq e^{- \frac 12 \gamma \|y-\aa_1(t)\|} e^{-\frac 14 \|\aa\|}.$$
Thus,
\bee
&& 
\int_{\|y-z\|\leq \frac 12 \|\aa\|} e^{- \gamma \|y-\aa_1(t)\|} |\e(t,y)| e^{- \gamma \|z-\aa_2(t)\|} |\e(t,z)|\frac 1{\|y-z\|} dy dz\\
&& \leq C e^{-\frac 14 \|\aa\|} \int_{\|y-z\|\leq \frac 12 \|\aa\|} e^{- \frac 12 \gamma \|y-\aa_1(t)\|} |\e(t,y)| 
|\e(t,z)| dy dz
\leq C \|\aa\|^{\frac 12} e^{-\frac 14 \|\aa\|} \|\e\|_{L^2}^2.
\eee

Gathering these estimates, we obtain
$$
{\cal{G}}(\e) \geq  \sum_{j=1}^2 
{\cal{G}}_{j,\psi}(\e) - \frac C{t^{\frac 23}} \|\e\|_{H^1}^2.
$$
Using (\ref{eeb}) and $\psi_1+\psi_2=1$, we find
${\cal{G}}(\e) \geq c \|\e\|_{H^1}^2 -\frac C{t^{\frac 23}} \|\e\|_{H^1}^2$
which gives the desired result for $T_0$ large enough.


\end{document}